\documentclass[12pt,a4paper]{article}
\usepackage{graphicx} 
\usepackage{epsfig}
\usepackage{amssymb}
\usepackage{amsmath}
\usepackage{amsfonts}
\usepackage{mathrsfs}
\usepackage[english]{babel}   
\usepackage{fancyhdr}
\usepackage{lastpage}
\usepackage[normal]{subfigure}
\usepackage{color}
\usepackage[dvipsnames]{xcolor}
\usepackage{inputenc}
\usepackage{dashrule}
\usepackage{float}

\usepackage{cancel}
\usepackage[normalem]{ulem}
\usepackage{breqn}

\newtheorem{comentario}{{\bf Note}}[section]

\def\QH{  { {\mathcal U}_{\mathcal Z}^H}}

\def\vu{\boldsymbol{u}}

\def\dxi{\Delta x_{i}}
\def\dx{\Delta x}
\def\dt{\Delta t}
\def\dxip{\Delta x_{i+\frac12}}

\def\a{\alpha}
\def\b{\beta}
\def\r{\rho}
\def\O{\Omega}
\def\D{\Delta}

\def\p{\partial}
\def\wtd{\widetilde}
\def\dsum{\displaystyle\sum}

\newcommand{\abs}[1]{\left\vert #1 \right\vert}

\definecolor{deepskyblue4}{rgb}{0,0.41,0.55}
\definecolor{dodgerblue4}{rgb}{0.06,0.31,0.55}
\definecolor{cadetblue2}{rgb}{0.56,0.9,0.93}
\definecolor{paleturquoise2}{rgb}{0.68,0.93,0.93}
\definecolor{deepskyblue3}{rgb}{0,0.6,0.8}
\definecolor{deepskyblue2}{rgb}{0,0.7,0.93}
\definecolor{royalblue3}{rgb}{0.23,0.37,0.8}
\definecolor{chartreuse3}{rgb}{0.4,0.8,0}
\definecolor{green3}{rgb}{0,0.8,0}
\definecolor{cyan3}{rgb}{0,0.8,0.8}
\definecolor{orange}{rgb}{1,0.65,0}
\definecolor{darkorange}{rgb}{1,0.55,0}
\definecolor{sienna1}{rgb}{1,0.51,0.28}
\definecolor{red2}{rgb}{0.93,0,0}
\definecolor{red1}{rgb}{1,0,0}
\definecolor{firebrick2}{rgb}{0.93,0.17,0.17}
\definecolor{yellow2}{rgb}{0.93,0.93,0}
\definecolor{gold1}{rgb}{1,0.84,0}
\definecolor{orchid}{rgb}{0.85,0.44,0.84}

\title{Flexible and efficient discretizations of multilayer models with variable density}

\author{Luca Bonaventura $^{(1)}$ and Jos\'e Garres-D\'iaz $^{(2)}$}

\begin{document}
\maketitle

\begin{center} 
{\small
$^{(1)}$ MOX -- Modelling and Scientific Computing, \\
Dipartimento di Matematica, Politecnico di Milano \\
Via Bonardi 9, 20133 Milano, Italy\\
{\tt luca.bonaventura@polimi.it}
}
\end{center}

\begin{center}
{\small $^{(2)}$
Departamento de Matem\'aticas, Universidad de C\'ordoba \\
 Campus de Rabanales, 14014, C\'ordoba, Spain\\
{\tt jgarres@uco.es}
}
\end{center}
\date{}


\vspace*{0.5cm}


\vspace*{0.5cm}

\pagebreak

\begin{abstract}
We show that the semi-implicit time discretization  approaches previously introduced
 for multilayer shallow water models for the barotropic case
 can be also applied to the variable density case with Boussinesq approximation. Furthermore, also for the variable density equations,
  a variable number of layers can  be used, so as to achieve greater flexibility and efficiency of the resulting multilayer approach.
  An analysis of the linearized system, which allows to derive linear stability parameters in simple configurations, and the resulting spatially semi-discretized equations are presented.
 A number of numerical experiments demonstrate the effectiveness of the proposed approach.
\end{abstract}

\section{Introduction}
\label{intro}

\indent

Multilayer shallow water models  \cite{audusse:2005,audusse:2008,audusse:2011,audusse:2011b}
have become quite popular in the last two decades to reduce the computational cost
of river and  coastal flow simulations. A version
of these models was derived in \cite{fernandez:2014}  from the   full Navier-Stokes system, by assuming a discontinuous profile of velocity, showing that the solution of the multilayer model is  a particular weak solution of the full Navier-Stokes system. 
In our previous work \cite{bonaventura:2018b}, we have shown that, in the barotropic, constant density and hydrostatic case, these models can be made more computationally
efficient by two complementary strategies. On the one hand, a classical semi-implicit time discretization can be employed
to remove the time step restriction based on the external gravity wave celerity, which adversely affects efficiency in low Froude
number regimes. On the other hand, we showed that multilayer shallow water models can also use  different numbers of layers
in different mesh locations, so as to reduce the computational cost and to allow for a more efficient allocation of the degrees of freedom as well as  for adaptive strategies.

Nevertheless, our previous work has an important limitation when simulating realistic flows, since the density is there assumed to be constant, as commented before. In variable density flows under the action of gravity, stratifications effects and internal gravity waves arise. Thus, the density field can have a strong effect over the dynamics of realistic coastal flows. In particular, it makes increase the vertical structure of the fluid because of the density variations (see e.g. \cite{csanady:1981,gill:1982}). Therefore the use of multilayer models is relevant in this context to account for this vertical structure. Several multilayer models and numerical strategies have already been proposed to deal with variable density flows. In \cite{fernandez:2013,morales:2017,burger:2019}, multilayer models with variable density  due to suspended sediments are introduced. In \cite{audusse:2011b}, the multilayer system \cite{audusse:2011} is extended to the variable density case, where the density is constant in each layer but may be different across the layers. This model is solved using a kinetic scheme, using again an explicit time discretization. Recently, a robust second-order explicit scheme for the multilayer systems with variable density, deduced from \cite{fernandez:2014}, has been proposed in \cite{guerrero:2020}. 

In the present paper, we extend the findings of \cite{bonaventura:2018b} to the hydrostatic, variable density case in the  Boussinesq regime.	This works is different from the results presented in  \cite{audusse:2011b,guerrero:2020} in several main aspects. Firstly, these papers use  explicit discretizations, while we use a semi-implicit time discretization to make more efficient the multilayer method for variable density flows (for the first time to our knowledge), which allows   to   reduce  significantly the computational cost in realistic simulations in the subcritical regime. Furthermore,  in the present work the number of vertical layers is no longer constant, leading to a more flexible and efficient discretization. We also work making the Boussinesq approximation, so that  the resulting system and also its numerical approximation are conceptually simpler, while still being applicable to coastal flow modelling.
 In addition, the proposed model and that in \cite{guerrero:2020} differ from \cite{audusse:2011b} in the procedure to obtain the multilayer system. Concretely, the vertical velocity and the momentum transference terms are different, and the solutions of our system yield a particular weak solution of the full Navier-Stokes system. In the present work we also study the linearization of the proposed model equations,
 showing how the presence of complex eigenvalues and the resulting loss of hyperbolicity are strictly related to the nontrivial vertical structure.
%
	
	The paper is organized as follows: in section \ref{se:model}, we derive the equations defining the variable density multilayer shallow water models. 
  In section \ref{se:linear}, the corresponding linearized equation are derived. In section \ref{se:spatial}, the  spatial semi-discretization is introduced,
while the time semi-discretization approach is introduced in section \ref{se:time}.
   Results of a number of numerical experiments 
are reported in section \ref{se:tests}, showing the significant efficiency gains that can be achieved by
the proposed techniques.    Conclusions and perspectives for future work are 
presented in section \ref{se:conclu}. 
 
\bigskip

\section{Variable density multilayer system}
\label{se:model}
We start here from the same  multilayer system as in the previous work \cite{bonaventura:2018b} and for convenience we recall the multilayer notation introduced there. The computational domain is divided in $N$ shallow vertical layers $\O_\a,\ \a=1,\dots,N$, where the upper and lower interfaces of layer $\O_{\a}$ are $\Gamma_{\a+\frac{1}{2}}$ and $\Gamma_{\a-\frac{1}{2}}$. In particular, $\Gamma_{1/2} = b$ and $\Gamma_{N+\frac12} = \eta$ denote the topography and the free surface. As usual, $h_\a$ denotes the height of the layer $\a=1,\dots,N$, while $h$ is the total height of the fluid, i.e., $h=\sum_{\a=1}^{N}h_\a$. Thus, the interfaces $\Gamma_{\a+\frac{1}{2}}$ are written $z_{\a+\frac12} = b+\sum_{\b=1}^{\a}h_\a$ and the free surface level is $\eta = b+h$. Finally, given a function $f$ which is continuous at the interface $\Gamma_{\a+\frac12}$, its approximation there is denoted by $f_{\a+\frac12}$.

Once the notation is fixed, the multilayer system for a fluid with constant density $\rho_0\in\mathbb{R}^+$ is written as 
\begin{eqnarray}
\label{eq:FinalModel}
 \p_{t}h_{\a} &+& \p_{x}\left(h_{\a}u_{\a}\right) = G_{\a+\frac{1}{2}} - G_{\a-\frac{1}{2}}, \nonumber \\
 \p_{t}\left(h_{\a}u_{\a}\right) \;&+&\; \p_{x}\left( h_{\a}u_{\a}^2\right) 
\quad+\; gh_{\a}\p_{x} \left(b+h\right) \, =\, \dfrac{1}{\rho_0}\left(K_{\a-\frac{1}{2}} - K_{\a+\frac{1}{2}}\right) \nonumber \\[4mm]
\quad&+&\;\dfrac{1}{2} G_{\a+\frac{1}{2}}\left(u_{\a+1} + u_{\a}\right) \;-\; \dfrac{1}{2} G_{\a-\frac{1}{2}}\left(u_{\a} + u_{\a-1}\right),
\end{eqnarray}
for $\a=1,\dots,N$, where $\left(u_\a,w_\a\right)\in \mathbb{R}^2$ is the two-dimensional velocity in the layer $\a$, $K_{\a+\frac12}$ account for the stresses between the layers, and $G_{\a+\frac12}$ is the mass transfer term between the layers $\a$ and $\a+1$. These terms at the interfaces are defined through the jump condition at $\Gamma_{\a+\frac12}:$ 
\begin{equation}
\label{eq:mass_transf_ini}
G_{\a+\frac{1}{2}} = \p_{t}z_{\a+\frac{1}{2}} + u_{\a+1}\p_{x} z_{\a+\frac{1}{2}} - w_{\a+\frac{1}{2}}^{+} = \p_{t}z_{\a+\frac{1}{2}} + u_{\a}\p_{x}z_{\a+\frac{1}{2}} - w_{\a+\frac{1}{2}}^{-},
\end{equation}
see \cite{fernandez:2014} for details.
By defining the positive coefficient $l_\a$, $\a=1,\dots,N$ such that $\sum_{\a=1}^N l_\a= 1$ and $h_\a = l_\a h$, and using the mass conservation equation in system \eqref{eq:FinalModel}, the mass transfer  terms are written as
\begin{equation}\label{eq:cont_mass_G}
G_{\a+\frac{1}{2}} = \  \dsum_{\b=1}^\a \Biggl[ \p_{x}\left(hl_{\b}u_{\b}\right) -  l_{\b}\dsum_{\gamma=1}^N\p_{x}\left(l_{\gamma}h u_{\gamma}\right)\Biggr].
\end{equation}
Here we assume that there is no mass transfer  at the bottom and the free surface, i.e., $G_{\frac12} = G_{N+\frac12} =0$. Then, using that $\p_t b =0, $ the multilayer system can be expressed as a system with $N+1$ equations and unknowns ($\eta,u_1,\dots,u_N$) 
\begin{eqnarray}
\label{eq:FinalModel_u}
\p_{t}\eta  &+& \p_{x}\Biggl(  h\dsum_{\b=1}^N l_{\b}u_{\b}\Biggr) = 0, \nonumber \\[2mm]
l_\a h\p_{t} u_{\a} &+&  l_\a h u_{\a} \p_x u_{\a}  \,  + gl_\a h \p_{x} \eta   \nonumber \\[4mm]
&=&\dfrac{1}{\rho_0}\left(K_{\a-\frac{1}{2}} - K_{\a+\frac{1}{2}}\right)\; 
+ G_{\a+\frac{1}{2}}\Delta \wtd{u}_{\a+\frac 12} +G_{\a-\frac{1}{2}}\Delta \wtd{u}_{\a-\frac 12},  
\end{eqnarray}
for $\a=1,\cdots,N$, where  $\Delta \wtd{u}_{\a+\frac 12}=( {u}_{\a+1}-{u}_{\a})/2.$ The viscous terms are
\begin{subequations}
	\label{eq:FinalModel2}
	\begin{alignat}{2}
	K_{\a+\frac{1}{2}} &=- \nu_{\a+\frac{1}{2}} \QH_{\a+\frac12} ,\end{alignat}
\end{subequations}
where $\nu$ denotes the kinematic viscosity and $\QH_{\a+\frac12}$ is an approximation of $\p_z u_{\a}$ at $\Gamma_{\a+\frac12}$. In principle, any model can be chosen to appropriately define $\nu_{\a+\frac{1}{2}}$. The cases $K_{\frac12}$ and $K_{N+\frac12}$ have to be defined by friction coefficients at the bottom and the free surface (wind stress).  Notice that a non-conservative formulation has been used, since the methods
proposed in this work are most appropriate for  subcritical flows. 

The extension of the previous multilayer system to the case of flows with variable density is performed assuming that the Boussinesq approximation is valid \cite{gill:1982}, i.e., that the density variations are so  small  that
their effects can be accounted for only in the  computation of the pressure gradient. 
To this aim, we define the density in layer $\a $ as
$\wtd{\r_\a} = \r_0 + \r_\a '$, where $\r_0$ is a constant reference density and $\r_\a' $ is assumed to be small.
More precisely, in typical geophysical applications one has
 $|\r_\a'/\wtd{\r_\a} |\leq O(10^{-2})$.  Under this hypothesis, the  density
perturbations  can be accounted for only in the pressure term of the momentum equation, which must be 
computed by vertical  integration in the  multilayer framework. 
After some straightforward algebra, we obtain

\begin{equation}\label{eq:pres_dens_var}
\begin{array}{ll}
\displaystyle \int_{\Gamma_{\a-1/2}}^{\Gamma_{\a+1/2}} \p_x p_\a(z) \ dz &= \r_0\,g\,h_\a\,\p_x\eta\, +\, g\,h_\a\,\p_x\left(\dsum_{\b=\a+1}^N \r_b' h_\b\right) \,+\,\\ &+ g\,\r_\a'\,h_\a\,\p_x \left(b+\dsum_{\b=1}^{\a-1}h_\b\right) + g\,\p_x\left(\r_\a'\,\dfrac{h_\a^2}{2}\right),
\end{array}
\end{equation}
where $p_\a(z)$ is the pressure term of the layer $\a$. Therefore, 
the momentum equation in the multilayer system can be written  as
\begin{eqnarray}
&&\r_0 l_\a h\p_{t} u_{\a} \;+\;  \r_0 l_\a h u_{\a}\p_x u_{\a}  \,  +  \r_0 g l_\a h  \p_{x} \eta \nonumber \\[4mm]
&&\quad +gl_\a h\p_x\left(\dfrac{\r_\a'h_\a}{2} + \dsum_{\b=\a+1}^N \r_\b' h_\b\right) + gl_\a h\r_\a' \partial_x \left(b+\dsum_{\b=1}^{\a-1}h_\b+\dfrac{h_\a}{2}\right) \nonumber\\[4mm]
&& \quad = K_{\a-\frac{1}{2}} - K_{\a+\frac{1}{2}}\; 
+ \r_0\, G_{\a+\frac{1}{2}}\Delta\tilde{u}_{\a+\frac 12} +\r_0\,G_{\a-\frac{1}{2}}\Delta\tilde{u}_{\a-\frac 12}.  
\end{eqnarray}
Notice that the pressure term  \eqref{eq:pres_dens_var} can also be rewritten as 
\begin{equation}\label{eq:pres_dens_var2}
\begin{array}{l}
\displaystyle \int_{\Gamma_{\a-1/2}}^{\Gamma_{\a+1/2}} \p_x p_\a(z) \ dz =\left(\r_0+\r_{\a}'\right)\,g\,h_\a\,\p_x\eta\, +\\[3mm] 
\qquad g\,h_\a\,\left(\dsum_{\b=\a+1}^N h_\b\p_x\r_{\b}'+\dfrac{h_\a\p_x\r_{\a}'}{2}\right) + gh_\a \dsum_{\b=\a+1}^N\left(\r_{\b}'-\r_{\a}'\right)\p_x h_\b,
\end{array}
\end{equation}
which can be useful if  separate treatments are sought for  the density perturbation gradients and  the layer thickness gradients.
 In addition, an equation is required for the evolution of $\r_\a'$. The mass conservation equation for the total density reads
\begin{equation*}
\p_{t}\left(\wtd{\r_\a} h_{\a}\right) + \p_{x}(\wtd{\r_\a} h_{\a}u_{\a}) = \r_0\,G_{\a+\frac{1}{2}} - \r_0\,G_{\a-\frac{1}{2}} + G'_{\a+\frac{1}{2}} - G'_{\a-\frac{1}{2}},
\end{equation*}  
where
$$G'_{\a+\frac{1}{2}} = \dfrac{\r_{\a+1}'+\r_\a'}{2}\p_t z_{\a+1/2} + \dfrac{\r_{\a+1}'u_{\a+1}+\r_\a'u_\a}{2}\p_x z_{\a+1/2}- \dfrac{\r_{\a+1}'w_{\a+1/2}^+ +\r_\a'w_{\a+1/2}^-}{2}.$$
Finally, by using the mass equation in \eqref{eq:FinalModel}, we obtain 
\begin{equation}\label{eq:dens_var}
\p_{t}\left(\r_\a' h_{\a}\right) + \p_{x}(\r_\a' h_{\a}u_{\a}) = G'_{\a+\frac{1}{2}} - G'_{\a-\frac{1}{2}}.
\end{equation}  
Since one has
$
G_{\a+\frac12}^{'\,+} = \r_{\a+1}' G_{\a+\frac12}; $  $  G_{\a+\frac12}^{'\,-} = \rho_{\a}' G_{\a+\frac12} $ and
$ G_{\a+\frac12}' = G_{\a+\frac12}^{'\,-} = G_{\a+\frac12}^{'\,+},$
(see \eqref{eq:mass_transf_ini})
it is easy to verify that $$G'_{\a+\frac12}=\r'_{\a+\frac12}G_{\a+\frac12}, \quad\text{with}\quad \r'_{\a+\frac12} = (\r'_{\a+1}+\r'_{\a})/2.$$
Therefore, system \eqref{eq:FinalModel_u}-\eqref{eq:FinalModel2} is finally re-written as
\begin{eqnarray}
\label{eq:sistema_u_dens}
\p_{t}\eta  &+& \p_{x}\Biggl(  h\dsum_{\b=1}^N l_{\b}u_{\b}\Biggr) = 0,  \nonumber  \\[4mm]
  l_\a h\p_{t} u_{\a} \;&+&\;  l_\a h u_{\a}\p_x u_{\a}  \,  +   g l_\a h  \p_{x} \eta \nonumber \\[2mm]
\quad &+&gl_\a h\p_x\left(\dfrac{\r_\a h_\a}{2} + \dsum_{\b=\a+1}^N \r_\b h_\b\right) + gl_\a h\r_\a \partial_x \left(b+\dsum_{\b=1}^{\a-1}h_\b+\dfrac{h_\a}{2}\right) \nonumber \\[2mm]
 \quad &=&\dfrac{1}{\r_0}\left(K_{\a-\frac{1}{2}} - K_{\a+\frac{1}{2}}\right)\; 
+ G_{\a+\frac{1}{2}}\Delta\tilde{u}_{\a+\frac 12} +G_{\a-\frac{1}{2}}\Delta\tilde{u}_{\a-\frac 12} \nonumber \\[4mm]
 \p_{t}\left(\r_\a l_\a h\right) &+& \p_{x}(\r_\a l_\a h u_{\a}) = \r_{\a+1/2}G_{\a+\frac{1}{2}} - \r_{\a-1/2}G_{\a-\frac{1}{2}},    
\end{eqnarray}
for $\a=1,\cdots,N, $ where we have redefined for simplicity the perturbation density as $\r_\a=\r_\a'/\r_0.$ 
If the reformulated pressure gradient \eqref{eq:pres_dens_var2} is employed, the corresponding equations 
\begin{eqnarray}
\label{eq:sistema_u_dens2}
 \p_{t}\eta  &+& \p_{x}\Biggl(  h\dsum_{\b=1}^N l_{\b}u_{\b}\Biggr) = 0, \nonumber  \\[4mm]
  l_\a h\p_{t} u_{\a} \;&+&\;  l_\a h u_{\a}\p_x u_{\a}  \,  +   g l_\a h (1+ \r_{\a}) \p_{x} \eta  \nonumber \\[2mm]
  \qquad &+&g\,l_\a h\,\left(\dsum_{\b=\a+1}^N h_\b\p_x\r_{\b}+\dfrac{h_\a\p_x\r_{\a}}{2}\right) + g l_\a h \dsum_{\b=\a+1}^N\left(\r_{\b}-\r_{\a}\right)\p_x h_\b \nonumber\\[2mm]
 \quad &=&\dfrac{1}{\r_0}\left(K_{\a-\frac{1}{2}} - K_{\a+\frac{1}{2}}\right)\; 
+ G_{\a+\frac{1}{2}}\Delta\tilde{u}_{\a+\frac 12} +G_{\a-\frac{1}{2}}\Delta\tilde{u}_{\a-\frac 12}, \nonumber \\[4mm]
 \p_{t}\left(\r_\a l_\a h\right) &+& \p_{x}(\r_\a l_\a h u_{\a}) = \r_{\a+1/2}G_{\a+\frac{1}{2}} - \r_{\a-1/2}G_{\a-\frac{1}{2}},    
\end{eqnarray}
are obtained for $\a=1,\cdots,N.$ 
 Notice that, as discussed \cite{fernandez:2014}, the above introduced equations can be regarded as a vertical discretization
 of the hydrostatic Navier-Stokes equations. This inevitably leads to the possibility that the  proposed equations
 fail to constitute an hyperbolic system, since the hydrostatic (or primitive) Navier-Stokes equations are well known
 to lose hyperbolicity as a consequence of the hydrostatic approximation, see e.g. the classical analysis in \cite{oliger:1978}.

\section{Linear analysis}
\label{se:linear}

 We will now derive explicitly the linearization of equations  \eqref{eq:sistema_u_dens},
 in order to study the hyperbolicity of these equations at least in the linear regime and to
  carry out stability analyses and discuss the  efficiency of time discretizations
in the variable density case.
We consider for simplicity the inviscid case $\nu_{\a+\frac{1}{2}}=0  $ with constant number of layers $N $
across the computational mesh. We assume
that $h=H+h', $ $b=0 $ (which implies $\eta=h $),   $u_{\a}=U_{\a}+u'_{\a}, $
where $H,U_{\a} $ are constants. Concerning the density variables, on the one hand, coherently  with the 
Boussinesq approximation, these are already small perturbations of the reference density $\r_0. $
On the other hand, however, the impact of stratification on the propagation of linear waves
can be significant. Therefore, we will consider first the general
case  $\r_\a =\varrho_\a + \r'_\a, $ allowing for different reference densities in each layer and assuming
the products   $h' \varrho_\a, $ $u'_{\a} \varrho_\a $  to be non negligible. We will then consider the simplified case
$\varrho_\a=0, $ in which possible stratification effects are disregarded. 

Considering then the former situation and  
 disregarding terms of second order in the perturbations, the following linearized equations are obtained,
where the primes denoting the perturbation variables have been dropped for convenience:

\begin{eqnarray}
\label{eq:sistema_u_dens_lin}
&&\p_{t}h  + H\dsum_{\b=1}^N l_{\b}\p_{x}u_{\b} +\dsum_{\b=1}^N l_{\b}U_{\b} \p_{x}  h   = 0,  \nonumber \\[4mm]
&&  \p_{t} u_{\a} \;+\;    U_{\a}\p_x u_{\a}  \,  +   g\left [1+\varrho_\a  + \dsum_{\b=\a+1}^N l_\b (\varrho_\b -\varrho_\a) \right ]   \p_{x} h  
\nonumber \\
&& \ \ \ \ \ \ +g\dfrac{l_\a H}2   \p_x\r_\a  + gH \dsum_{\b=\a+1}^N l_\b \p_x\r_\b    \nonumber \\[2mm]
&&= \bar G_{\a+\frac{1}{2}}\frac{ {U}_{\a+1} -{U}_{\a} }{2Hl_\a} +\bar G_{\a-\frac{1}{2}}
\frac{{U}_{\a} -{U}_{\a-1}}{2Hl_\a}, \nonumber \\[4mm]
&& \p_{t} \r_\a  +   U_{\a}\p_{x}\r_\a   = \dfrac{\varrho_{\a+1}-\varrho_{\a}}{2Hl_\a}\bar G_{\a+\frac{1}{2}} + \dfrac{\varrho_{\a}-\varrho_{\a-1}}{2Hl_\a}\bar G_{\a-\frac{1}{2}},    
\end{eqnarray}
where now
\begin{equation}\label{eq:cont_mass_G_lin}
\bar G_{\a+\frac{1}{2}} = \  \dsum_{\b=1}^\a \Biggl[
l_{\b}H\p_{x} u_{\b}   +l_{\b}U_{\b}\p_{x}h -  l_{\b}\dsum_{\gamma=1}^Nl_{\gamma}\left(H \p_{x}u_{\gamma}+ U_{\gamma}\p_{x}h\right)\Biggr]
\end{equation}
for $\a=1,\dots,N-1$ and $\bar G_{\frac{1}{2}} =\bar G_{N+\frac{1}{2}} =0 $ as in the nonlinear case.
It is important to remark here that the same  equations \eqref{eq:sistema_u_dens_lin} also arise  
from linearization of the modified system \eqref{eq:sistema_u_dens2}.
For compactness, we now set, again for $\a=1,\dots,N-1,$

\begin{eqnarray}
&& \bar U=\dsum_{\gamma=1}^N l_{\gamma}U_{\gamma} \ \ \  \ \delta U_\b = U_\b - \bar U \ \ \ \
\ \ \overline{\delta U}_{\a}=\dsum_{\gamma=1}^\a  l_{\gamma} \delta U_{\gamma},   \nonumber \\
&& \delta_{\b,\a+\frac 12}U=\frac{ {U}_{\a+1} -{U}_{\a} }{2Hl_\b},  \ \ \ r_\a= \varrho_\a  + \dsum_{\b=\a+1}^N l_\b (\varrho_\b -\varrho_\a) \ \ \ \
 \end{eqnarray}
and we also define the matrix
 \begin{eqnarray}
M_{\alpha,\gamma} &=&  \quad l_{\gamma}\left(1-\dsum_{\b=1}^\a l_\b\right) \ \ \ \ \ \ \ {\rm for} \ \ \gamma \leq \alpha \nonumber \\
M_{\alpha,\gamma} &=& \quad -l_{\gamma}\dsum_{\b=1}^\a l_\b \ \ \ \ \ \ \ \ \ {\rm for} \ \ \gamma > \alpha,
 \end{eqnarray}
 so that
\begin{equation}\label{eq:cont_mass_G_lin2}
\bar G_{\a+\frac{1}{2}} =  \overline{\delta U}_{\a}\p_{x}h +H\dsum_{\gamma=1}^N M_{\alpha,\gamma} \p_{x} u_{\gamma}
\end{equation}
and the momentum equations can be rewritten as
\begin{eqnarray}
\label{eq:momeq_re}
&&  \p_{t} u_{\a} \;+\;    U_{\a}\p_x u_{\a}  \,  +   g(1+r_\a) \p_{x} h   +g\dfrac{l_\a H}2   \p_x\r_\a  + gH \dsum_{\b=\a+1}^N l_\b \p_x\r_\b    \nonumber \\[2mm]
&&= \left [ (\delta_{\a,\a+\frac 12}U)\overline{\delta U}_{\a} + ( \delta_{\a,\a-\frac 12}U)\overline{\delta U}_{\a-1}\right]   \p_{x}h \nonumber \\
&&+H    \dsum_{\gamma=1}^N \left [    (\delta_{\a,\a+\frac 12}U) M_{\alpha,\gamma}  + ( \delta_{\a,\a-\frac 12}U)  M_{\alpha-1,\gamma} \right] \p_{x} u_{\gamma}\nonumber \\
&&= -v_\a   \p_{x}h -  H \dsum_{\gamma=1}^N W_{\alpha,\gamma}\p_{x} u_{\gamma},
\end{eqnarray}
while the density equation reads

\begin{eqnarray}
\label{eq:denseq_re}
  \p_{t} \r_\a  +   U_{\a}\p_{x}\r_\a  
&=& \left [ (\delta_{\a,\a+\frac 12}\varrho)\overline{\delta U}_{\a}  +(\delta_{\a,\a-\frac 12}\varrho)\overline{\delta U}_{\a-1}\right]   \p_{x}h \nonumber \\
&+& H    \dsum_{\gamma=1}^N \left [    (\delta_{\a,\a+\frac 12}\varrho) M_{\alpha,\gamma}   +(\delta_{\a,\a-\frac 12}\varrho)  M_{\alpha-1,\gamma} \right] \p_{x} u_{\gamma}\nonumber \\
&=& -v^{\r }_{\a}   \p_{x}h -    H \dsum_{\gamma=1}^N W^{\r}_{\alpha,\gamma}\p_{x} u_{\gamma}.
\end{eqnarray}
 Setting now  $$\mathbf{q}=(h,\mathbf{u},\boldsymbol{\rho})=(h,u_1,\dots,u_N,\r_1,\dots,\r_N),    $$
the previous equations can be written as
$$
\p_{t} \mathbf{q}+  \mathbf{A}\p_{x} \mathbf{q}= \mathbf{0}, 
$$
where we have now defined
\begin{equation}
\mathbf{A}=\left[
\begin{array}{ccccccc}
\bar U &  l_{1}H  &\dots & l_{N}H&  0&\dots &0\\
gs_1&  U_1+HW_{1,1} & \dots & HW_{1,N}& gl_1 H/2 &  \dots&gl_N H \\
gs_2&  HW_{2,1} & \dots &HW_{2,N}&0&   \dots&gl_N H \\
\vdots&\vdots & \ddots& \vdots& \vdots&  \ddots&\vdots \\
gs_N& HW_{N,1}   &\dots & U_N+HW_{N,N} & 0 &\dots  &  g l_N H/2\\
v^{\r }_1 &  HW^\r_{1,1}   &\dots & HW^\r_{1,N}  & U_{1} & \dots& 0\\
v^{\r }_2 & HW^\r_{2,1}   &\dots & HW^\r_{2,N} & 0 & \dots& 0\\
\vdots&\vdots  & \ddots& \vdots& \vdots & \ddots&\vdots \\
v^{\r }_N &HW^\r_{N,1} & \dots &  HW^\r_{N,N}& 0& \dots& U_{N}\\
\end{array}
\right], \nonumber
\end{equation}
where now $s_\a=(1+r_\a +v_\a). $
More compactly, this can also be rewritten as
\begin{equation}
\mathbf{A}=\left[
\begin{array}{ccc}
\bar U &  H {\bf l}^T&  {\bf 0}\\
g{\bf s} & {\bf U}+H{\bf W} & gH{\bf T}\mathbf{L}\\
{\bf v}^{\r} &   H{\bf W}^{\r} & {\bf U} \\
\end{array}
\right], \nonumber
\end{equation}
where ${\bf l}=[l_1,\dots,l_N]^T, $ ${\bf 0}=[0,\dots,0]^T, $  ${\bf s} =[s_1,\dots,s_N]^T, $ $ {\bf v}^{\r}=[ v^{\r }_1,\dots, v^{\r }_N]^T. $
Note that ${\bf W}$ and ${\bf W}^{\r}$ are the contributions corresponding to mass transference terms multiplied by the differences of reference velocities and densities, respectively. Furthermore, 
 ${\bf L} $ denotes the diagonal matrix with elements $l_i $ on the main diagonal, 
${\bf U} $ denotes the diagonal matrix with elements $U_a $ on the main diagonal
 and $\mathbf{T}$ denotes the upper triangular matrix such that

\begin{equation}
\mathbf{T}=\left[
\begin{array}{cccc}
 1/2 & 1               &  \dots  & 1 \\
  0                       &    1/2 &  \dots  & 1 \\
         \vdots& \vdots & \ddots&\vdots \\
     0 &0 &\dots  &   1/2\\
\end{array}
\right] {\rm so \ that \ }
\mathbf{T}\mathbf{L}=\left[
\begin{array}{cccc}
 l_1/2 & l_2                 &  \dots  & l_N  \\
  0                       &    l_2 /2 &  \dots  & l_N \\
         \vdots& \vdots & \ddots&\vdots \\
     0 &0 &\dots  &   l_N/2\\
\end{array}
\right]. \nonumber
\end{equation}
The structure of  ${\bf A} $ is simpler if special cases are considered. For example,
it can be immediately observed that, for a constant reference velocity profile, e.g. $U_\a=\bar U =U, $
one has  $v_\a=  0 $ and ${\bf W}={\bf 0}.$ If the reference density values are also
taken to be zero, one then has $v^{\r }_\a=0, $ ${\bf W}^{\r}={\bf 0} $ and
${\bf s} ={\bf e}=[1,\dots,1]^T, $ thus yielding

\begin{equation}
\mathbf{A}=\left[
\begin{array}{ccc}
  U &  H {\bf l}^T&  {\bf 0}\\
g{\bf e} & {\bf U} & gH{\bf T}\mathbf{L}\\
{\bf 0} &  {\bf 0} & {\bf U} \\
\end{array}
\right]. \nonumber
\end{equation}
In this particular case, matrix $\mathbf{A} $ has eigenvalues $U\pm\sqrt{gH} $ and $U $ (with multiplicity $2N-1$),
independently of the layer distribution, so that the system is hyperbolic. If instead the  reference velocity
and density  profiles are both constant, but ${\r }_\a={\r }\neq 0, $ one still has has $v^{\r }_\a=0, $ ${\bf W}^{\r}={\bf 0}, $
so that
\begin{equation}
\mathbf{A}=\left[
\begin{array}{ccc}
  U &  H {\bf l}^T&  {\bf 0}\\
g{\bf s} & {\bf U} & gH{\bf T}\mathbf{L}\\
{\bf 0} &  {\bf 0} & {\bf U} \\
\end{array}
\right] \nonumber
\end{equation}
with $s_\a=1+\rho, \ \a=1,\dots, N, $ so that the matrix $\mathbf{A} $ has eigenvalues $U\pm\sqrt{g(1+\rho)H} $ and $U $ (with multiplicity $2N-1$),
independently of the size of the layers. On the other hand, if ${\r }_\a=   0 $ but the reference velocity profile is not constant,
it follows

\begin{equation}
\mathbf{A}=\left[
\begin{array}{ccc}
\bar U &  H {\bf l}^T&  {\bf 0}\\
g{\bf e} & {\bf U}+H{\bf W} & gH{\bf T}\mathbf{L}\\
{\bf 0}  &   {\bf 0}  & {\bf U} \\
\end{array}
\right]. \nonumber
\end{equation}
The characteristic polynomial of $\mathbf{A}$ reads then
\begin{equation}
{\rm det}(\mathbf{A}-\lambda\mathbf{I})= {\rm det}(\mathbf{U}-\lambda\mathbf{I})
{\rm det}\left(\left[
\begin{array}{cc}
\bar U -\lambda &  H {\bf l}^T \\
g{\bf e} & {\bf U}+H{\bf W}- \lambda\mathbf{I} \\ \end{array}
\right]\right).
\end{equation}
By application of Banachiewicz-like decompositions, see e.g. \cite{silvester:2000,tian:2009},
it follows that the determinant of the last matrix is equal to
$$
(\bar U -\lambda)   {\rm det}\left({\bf U}+H{\bf W}- \frac{\lambda-gH}{\bar U -\lambda}    \mathbf{I}\right).
$$
This implies that, if $z $ denotes an eigenvalue of $ {\bf U}+H{\bf W},$
the eigenvalues $\lambda $ of  $\mathbf{A}$ must satisfy the equation
$$
(z-\lambda)(\bar U -\lambda)-gH=0.
$$
Therefore, is it sufficient that the non symmetric part of ${\bf W} $ is large enough to yield
complex eigenvalues for  $\mathbf{A}$ as well. As remarked in the introduction,
this should not be regarded as a deficiency of the multilayer model, but rather as a consequence
of its being a convergent approximation of a three-dimensional hydrostatic flow. Let us remark that this situation is related to large deviations from the constant velocity profile, and therefore to large values of the mass transference terms $G_{\a+1/2}$. To our experience, moderate values are obtained in simulations in the hydrostatic regime. Therefore, in practice the linearized system \eqref{eq:sistema_u_dens_lin} is expected to be hyperbolic.
\section{Spatial discretization}
\label{se:spatial}

We now consider a spatial discretization for system \eqref{eq:sistema_u_dens}, which extends to the variable density
case the finite volume approach presented in  \cite{bonaventura:2018b}. We only recall the main features of the 
discretization, referring to our previous paper for the full description and   focusing on the novel terms coming from the variable density pressure and the evolution of the perturbations of density in each layer. 
It should be remarked that many other options can also be considered, such as finite difference or finite element methods.

\begin{figure}[!h]
	\begin{center}
		\includegraphics[width=0.6\textwidth]{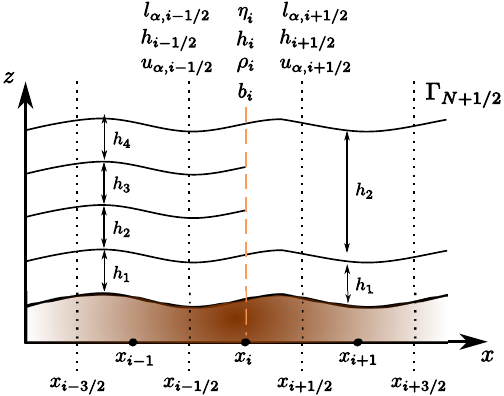}
	\end{center}
	\caption{\label{fig:Multilayers_part_den} \footnotesize\it{Sketch of the domain and of its subdivision in a variable number of   layers.}}
\end{figure}

We consider the usual finite volume description of the horizontal domain, which is subdivided into control volumes $V_i=(x_{i-1/2},x_{i+1/2})$, with centers  $x_i = (x_{i+1/2}+x_{i-1/2})/2$, for $i=1,\dots,M$. The distances between $x_i$ and $x_{i+1}$ is denoted $\dx_{i+1/2}$ and the length of the control volume is $\dxi = x_{i+1/2} - x_{i-1/2}$. A staggered mesh is considered, where the discrete free surface and density variables $\eta_i,\r_i$ are defined at the centers of the control volumes, $x_i$, while the discrete velocities $u_{\a,i+\frac12}$ are defined at the interfaces, $x_{i+1/2}$ (see figure \ref{fig:Multilayers_part_den}).

As discussed in \cite{bonaventura:2018b},  the positive coefficients $l_{\a,i+\frac12}$ are also defined
at half-integer locations, so that the number of layers $N_{i+\frac12} $ is also specified at these locations and may
vary across the discrete mesh. For simplicity, the vertical mesh is assumed to be conformal, i.e., either a layer splits in several of some layers merge into a single one, and only a transition between cells with different number of layers is allowed in a 3-point stencil, i.e., it is not possible to have two consecutive transitions. These two hypothesis reduce the complexity of the implementation, namely in the case of advection terms in the momentum equation.  
The number of layers at integer locations is defined as $N_i=\max\{N_{i-\frac12} ,N_{i+\frac12}\} $
and the discrete layer thickness coefficients at integer locations  $l_{\a,i} $ are taken to be equal to those at the neighbouring half-integer location with larger number of layers, that is, assuming for example $N_{i-\frac12} < N_{i+\frac12}, $ we take
$l_{\a,i} =l_{\a,i+\frac12}.$

Given these definitions, the spatial semi-discretization reads then

\begin{equation}
 \dfrac{d \eta_i}{dt} = -\dfrac{1}{\dxi} \dsum_{\b=1}^N \Big( l_{\b,i+\frac12}h_{i+\frac12}u_{\b,i+\frac12} \,-\, l_{\b,i-\frac12}h_{i-\frac12}u_{\b,i-\frac12} \Big)
 \end{equation}
 \begin{eqnarray}
 \dfrac{d u_{\a,i+\frac12}}{dt} &= &-\dfrac{1}{\dxip} l_{\a,i+\frac12}h_{i+\frac12}\,g\left(\eta_{i+1}-\eta_i\right)\nonumber\\[4mm]
\qquad &+& \bigg(\nu_{\a+\frac12,i+\frac12}\,\dfrac{u_{\a+1,i+\frac12}-u_{\a,i+\frac12}}{l_{\a+\frac12,i+\frac12}h_{i+\frac12}} - \nu_{\a-\frac12,i+\frac12}\,\dfrac{u_{\a,i+\frac12}-u_{\a-1,i+\frac12}}{l_{\a-\frac12,i+\frac12}h_{i+\frac12}} \bigg)\nonumber \\
 &-&l_{\a,i+\frac12}\,h_{i+\frac12}\,\mathcal{A}_{\a,i+\frac12}^{u}+ \dfrac{1}{\dxip}\Big(\Delta \wtd{u}_{\a+\frac12,i+\frac12}
\mathcal{G}_{\a+\frac12,i+\frac12}+\Delta \wtd{u}_{\a-\frac12,i+\frac12}\mathcal{G}_{\a-\frac12,i+\frac12}\Big) \nonumber\\[4mm]
&-&\, \dfrac{l_{\a,i+\frac12}h_{i+\frac12}}{\dx_{i+1/2}} g \left(
\dsum_{\b=\a+1}^N l_{\b,i+\frac12}\left(\r_{\b,i+1}h_{i+1}-\r_{\b,i}h_{i}\right) + l_{\a,i+\frac12}\dfrac{\r_{\a,i+1}h_{i+1} - \r_{\a,i}h_{i}}{2}
\right)\nonumber \\[4mm]
&-&\, \dfrac{l_{\a,i+\frac12}h_{i+\frac12}}{\dx_{i+1/2}} g \r_{\a,i+\frac12}\left(b_{i+1}-b_i +
\left(h_{i+1}-h_{i}\right)\left(\dsum_{\b=1}^{\a-1} l_{\b,i+\frac12} + \dfrac{l_{\a,i+\frac12}}{2}\right)
\right),
\end{eqnarray}
 \begin{eqnarray}
 \dfrac{d \rho_{\a,i}}{dt} &= &-\dfrac{1}{l_{\a,i}\dxip} \left(l_{\a,i+\frac12}h_{i+\frac12}\r_{\a,i+\frac12}u_{\a,i+\frac12} - l_{\a,i-\frac12}h_{i-\frac12}\r_{\a,i-\frac12}u_{\a,i-\frac12}\right)\nonumber \\
&+& \dfrac{1}{l_{\a,i}} \left(\r_{\a+\frac12,i}G_{\a+\frac12,i}-\r_{\a-\frac12,i}G_{\a-\frac12,i}\right)
\end{eqnarray}
for $\a=1,\dots,N. $ Here $h_{i+\frac12},\r_{\a,i+\frac12}, l_{\a,i}$ are the defined as the upwind values
and  $\mathcal{A}_{\a,i+\frac12}^{u}$ denotes the discretization of the momentum advection term, which is carried out
by a first or second order upwind method. More specifically, we use the same second-order upstream based second order finite difference approximation as in \cite{bonaventura:2018b}, i.e. 
	$
	\Big(u_\a\p_x u_\a\Big){\bigg|_{i+1/2}} = u_{\a,i+1/2}\, \p_x u_{\a}{\Big|_{i+1/2}}, $
	where
	$$
	  \p_x u_{\a}{\Big|_{i+1/2}} = \left\{\begin{array}{lcc}
	\dfrac{u_{\a,i-\frac32} - 4u_{\a,i-\frac12}+3u_{\a,i+\frac12}}{2\Delta x_{i+\frac12}} & \mbox{if} & u_{\a,i+\frac12}>0,\\[6mm]
	-\dfrac{u_{\a,i+\frac52} - 4u_{\a,i+\frac32}+3u_{\a,i+\frac12}}{2\Delta x_{i+\frac12}} & \mbox{if} & u_{\a,i+\frac12}<0.
	\end{array}\right.
	$$

\section{IMEX-ARK2 time discretization}
\label{se:time}

This spatial discretization may be coupled with any time discretization.
As in the previous work, we use an explicit third order Runge-Kutta method for the reference solutions and one of the semi-implicit scheme presented in \cite{bonaventura:2018b}, the second order IMplicit-EXplicit Additive Runge-Kutta method (IMEX-ARK2) including the terms coming to the variable density pressure terms. Moreover, we have an additional equation for the perturbations of the density, which must be discretized accordingly.

We extend here the method presented in \cite{bonaventura:2018b} to the variable density case. To introduce IMEX-type methods,
we write the ODE system as 
 $$\mathbf{y}' = \mathbf{f}_{s}(\mathbf{y},t) +\mathbf{f}_{ns}(\mathbf{y},t), $$
  where the $s$ and $ns$ subscripts denote the stiff
and non-stiff components of the system, respectively. In our case, we have  
$$
\begin{array}{ll}
\mathbf{f}_{s}^0 = &-\dfrac{1}{\dxi} \dsum_{\b=1}^N \Big( l_{\b,i+\frac12}h_{i+\frac12}u_{\b,i+\frac12} \,-\, l_{\b,i-\frac12}h_{i-\frac12}u_{\b,i-\frac12} \Big);\\[4mm]
\mathbf{f}_{s}^{2\a-1} = &-\dfrac{1}{\dxip} l_{\a,i+\frac12}h_{i+\frac12}\,g\left(\eta_{i+1}-\eta_i\right)\\[4mm]
\qquad &+ \bigg(\nu_{\a+\frac12,i+\frac12}\,\dfrac{u_{\a+1,i+\frac12}-u_{\a,i+\frac12}}{l_{\a+\frac12,i+\frac12}h_{i+\frac12}} - \nu_{\a-\frac12,i+\frac12}\,\dfrac{u_{\a,i+\frac12}-u_{\a-1,i+\frac12}}{l_{\a-\frac12,i+\frac12}h_{i+\frac12}} \bigg),\\[4mm]
\mathbf{f}_{s}^{2\a} = &-\dfrac{1}{l_{\a,i}\dxip} \left(l_{\a,i+\frac12}h_{i+\frac12}\r_{\a,i+\frac12}u_{\a,i+\frac12} - l_{\a,i-\frac12}h_{i-\frac12}\r_{\a,i-\frac12}u_{\a,i-\frac12}\right),
\end{array}
$$
and
$$
\begin{array}{l}
\mathbf{f}_{ns}^0 = 0;\\[4mm]
\mathbf{f}_{ns}^{2\a-1} = -l_{\a,i+\frac12}\,h_{i+\frac12}\,\mathcal{A}_{\a,i+\frac12}^{u}+ \dfrac{1}{\dxip}\Big(\Delta \wtd{u}_{\a+\frac12,i+\frac12}
\mathcal{G}_{\a+\frac12,i+\frac12}+\Delta \wtd{u}_{\a-\frac12,i+\frac12}\mathcal{G}_{\a-\frac12,i+\frac12}\Big) \\[4mm]
-\, \dfrac{l_{\a,i+\frac12}h_{i+\frac12}}{\dx_{i+1/2}} g \left(
\dsum_{\b=\a+1}^N l_{\b,i+\frac12}\left(\r_{\b,i+1}h_{i+1}-\r_{\b,i}h_{i}\right) + l_{\a,i+\frac12}\dfrac{\r_{\a,i+1}h_{i+1} - \r_{\a,i}h_{i}}{2}
\right)\\[4mm]
-\, \dfrac{l_{\a,i+\frac12}h_{i+\frac12}}{\dx_{i+1/2}} g \r_{\a,i+\frac12}\left(b_{i+1}-b_i +
\left(h_{i+1}-h_{i}\right)\left(\dsum_{\b=1}^{\a-1} l_{\b,i+\frac12} + \dfrac{l_{\a,i+\frac12}}{2}\right)
\right),\\[4mm]
\mathbf{f}_{ns}^{2\a} = \dfrac{1}{l_{\a,i}} \left(\r_{\a+\frac12,i}G_{\a+\frac12,i}-\r_{\a-\frac12,i}G_{\a-\frac12,i}\right),
\end{array}
$$
for $\a=1,\dots,N$, where $h_{i+\frac12},\r_{\a,i+\frac12}, l_{\a,i}$ are the upwind values. Note that the stiff and non-stiff parts for the extra equation are defined in order to be consistent with the continuity equation in the sense of \cite{gross:2002}. 

Then, the $s-$stage IMEX-ARK2 method can be defined as follows. For $l=1,2,3$: 
\begin{equation}
\label{imex-ark}
\begin{array}{rcl}
\mathbf{u}^{(l)}=\mathbf{u}^n&+&\Delta t  \dsum_{m=1}^{l-1} \bigg( a_{lm}\mathbf{f}_{ns}(\mathbf{u}^{(m)},t+c_m \Delta t) \\[4mm]
&+&\wtd{a}_{lm}\mathbf{f}_{s}(\mathbf{u}^{(m)},t+c_m\Delta t) \bigg) + \Delta t \, \wtd{a}_{ll} \, \mathbf{f}_{s}(\mathbf{u}^{(l)},t+c_l\Delta t).
\end{array}
\end{equation}
and the updates values $\vu^{n+1}$ are computed as
$$
\mathbf{u}^{n+1}=\mathbf{u}^n+\Delta t \sum_{l=1}^{3}b_{l}(\mathbf{f}_{ns}(\mathbf{u}^{(l)},t+c_l \Delta t)+\mathbf{f}_{s}(\mathbf{u}^{(l)},t+c_l\Delta t)).
$$
Coefficients $a_{lm}, \wtd{a}_{lm}, c_l$ and $b_l$  are given to obtain a consistent method satisfying specific order and coupling conditions.
We use the IMEX method proposed in \cite{giraldo:2013}, whose coefficients are in the Butcher tableaux, table \ref{ark2_butch_e} and \ref{ark2_butch_i} for the explicit and implicit method, respectively. 
The implicit part of this method matches with the TR-BDF2 scheme (see \cite{hosea:1996}), which is $L$-stable, and for the explicit part we have to respect the stability condition given by the Courant number $C_{vel}$ \eqref{eq:cfl_vel}. The coefficients of this part were proposed in
\cite{giraldo:2013}.
\begin{table}[!h]
	\begin{center}
		\begin{tabular}{c|ccc}
			0 & 0 & &  \\
			$2 \mp \sqrt{2}$ & $2 \mp \sqrt{2}$ & 0 & \\
			1 & $1-(3+2\sqrt{2})/6$ & $ (3+2\sqrt{2})/6$ & 0 \\
			\hline
			& $\pm \frac{1}{2\sqrt{2}}$ & $\pm \frac{1}{2\sqrt{2}}$ & $1 \mp \frac{1}{\sqrt{2}}$
		\end{tabular}
		\begin{tabular}{llll}
			& \multicolumn{1}{}{}  &     &  \\
			& \multicolumn{1}{l|}{}    &  &  \\
			& \multicolumn{1}{l|}{$c_l$} & $a_{lm}$ &    \\ \cline{2-4}
			& \multicolumn{1}{l|}{}  & $b_l$ &    \\
		\end{tabular}
	\end{center}
	\caption{\it Butcher tableaux of the explicit ARK2 method}
	\label{ark2_butch_e}
\end{table}

\begin{table}[h!]
	\begin{center}
		\begin{tabular}{c|ccc}
			0 & 0 & & \\
			$2 \mp \sqrt{2}$ &  $1 \mp \frac{1}{\sqrt{2}}$ &  $1 \mp \frac{1}{\sqrt{2}}$ & \\
			1 & $\pm \frac{1}{2\sqrt{2}}$ & $\pm \frac{1}{2\sqrt{2}}$ & 1 $\mp \frac{1}{\sqrt{2}}$ \\
			\hline
			& $\pm \frac{1}{2\sqrt{2}}$ & $\pm \frac{1}{2\sqrt{2}}$ & $1 \mp \frac{1}{\sqrt{2}}$
		\end{tabular}
		\begin{tabular}{llll}
			& \multicolumn{1}{}{}    &  &  \\
			& \multicolumn{1}{l|}{}    &  &  \\
			& \multicolumn{1}{l|}{$c_l$} & $\wtd{a}_{lm}$ &    \\ \cline{2-4}
			& \multicolumn{1}{l|}{}  & $b_l$ &    \\
		\end{tabular}
	\end{center}
	\caption{\it Butcher tableaux of the implicit ARK2 method}
	\label{ark2_butch_i}
\end{table}

Next, we detail the IMEX discretization for the equation of the evolution of the perturbations of the density. The mass and the momentum equations are discretized at each stage as detailed in \cite{bonaventura:2018b}, then we refer the reader to previous work for that discretized equations. 

In order to avoid solving an extra linear system to find $\r_{\a,i+\frac12}^{n,k},\ k=2,3$, we propose a modified IMEX discretization. It consist of linearize of the density in the second and third stages of the IMEX scheme, as we detail in the following.

For the first stage, we define $ \eta_{i}^{n,1}=\eta_{i}^{n}, $  $\r_{\a,i}^{n,1}=\r_{\a,i}^{n},$ and $u_{\a,i+\frac12}^{n,1}=u_{\a,i+\frac12}^{n},$ 
respectively. For the second stage, we have
\begin{equation*}\label{eq:evol_rho_imex1}
\begin{array}{l}
l_{\a,i}h_{i}^{n,2}\r_{\a,i}^{\,n,2} =   l_{\a,i}h_{i}^{n}\r_{\a,i}^{\,n}\\ [4mm]
\qquad- \dfrac{\dt}{\dxi}\left( l_{\a,i+\frac12} h^{n}_{i+\frac12} \r_{\a,i+\frac12}^{n,1} u^{*,2}_{\a,i+\frac12}
\,-\,l_{\a,i-\frac12} h_{i-\frac12}^{n} \r_{\a,i-\frac12}^{n,1}u^{*,2}_{\a,i-\frac12} \right)  \\[4mm]
\qquad+\,  a_{21}\dt\left(\r_{\a+\frac12,i}^{n,1} G_{\a+\frac12,i}^{n,1} - \r_{\a-\frac12,i}^{n,1} G_{\a-\frac12,i}^{n,1}\right)\, ,
\end{array}
\end{equation*}
where
$   u_\a^{*,2} =  \wtd{a}_{22} u_\a^{n,2} + \wtd{a}_{21}u_\a^{n,1}$. Here, $\r_{\a,i+\frac12}^{n,1}$ is computed as the upwind value depending on $u_\a^{*,2}$.
For the third stage, 
\begin{equation*}\label{eq:evol_rho_imex2}
\begin{array}{l}
l_{\a,i}h_{i}^{n,3}\r_{\a,i}^{\,n,3} = l_{\a,i}h_{i}^{n}\r_{\a,i}^{\,n}\\ [4mm]
\qquad-\, \dfrac{\dt}{\dxi}\left( l_{\a,i+\frac12} h^{n}_{i+\frac12} \r_{\a,i+\frac12}^{n,2} u^{*,3}_{\a,i+\frac12}
\,-\,l_{\a,i-\frac12} h_{i-\frac12}^{n} \r_{\a,i-\frac12}^{n,2}u^{*,3}_{\a,i-\frac12} \right) \\[4mm]
\qquad+\,  a_{32}\dt\left(\r_{\a+\frac12,i}^{n,2} G_{\a+\frac12,i}^{n,2} -  \r_{\a-\frac12,i}^{n,2} G_{\a-\frac12,i}^{n,2}\right)\\[4mm]
\qquad-\, \wtd{a}_{31}\dfrac{\dt}{\dxi}\left( l_{\a,i+\frac12} h^{n}_{i+\frac12} \r_{\a,i+\frac12}^{n,1} u^{n,1}_{\a,i+\frac12}
\,-\,l_{\a,i-\frac12} h_{i-\frac12}^{n} \r_{\a,i-\frac12}^{n,1}u^{n,1}_{\a,i-\frac12} \right) \\[4mm]
\qquad+\,  a_{31}\dt\left(\r_{\a+\frac12,i}^{n,1} G_{\a+\frac12,i}^{n,1} - \r_{\a-\frac12,i}^{n,1} G_{\a-\frac12,i}^{n,1}\right),
\end{array}
\end{equation*}
where now
$   u_\a^{*,3} =  \wtd{a}_{32} u_\a^{n,3} + \wtd{a}_{31}u_\a^{n,2}$, and $\r_{\a,i+\frac12}^{n,2}$ the upwind value depending on $u_\a^{*,3}$.  
Finally, the assembly of the solution at time level $n+1$ is  
\begin{equation}\label{eq:evol_rho_imex_final}
\begin{array}{l}
l_{\a,i}h_{i}^{n+1}\r_{\a,i}^{\,n+1} =  l_{\a,i}h_{i}^{n}\r_{\a,i}^{\,n} \\[4mm]
\quad- \,\dfrac{\dt}{\dxi} \dsum_{j=1}^{3} b_{j} \left( l_{\a,i+\frac12} h^{n}_{i+\frac12} \r_{\a,i+\frac12}^{n,j} u^{n,j}_{\a,i+\frac12} \right .
-\,\left. l_{\a,i-\frac12} h_{i-\frac12}^{n} \r_{\a,i-\frac12}^{n,j}u^{n,j}_{\a,i-\frac12} \right)\\[4mm]
\quad+\,   \dt  \dsum_{j=1}^{3} b_{j} \left(\r_{\a+\frac12,i}^{n,j} G_{\a+\frac12,i}^{n,j}
- \r_{\a-\frac12,i}^{n,j} G_{\a-\frac12,i}^{n,j}\right). 
\end{array}
\end{equation}
Notice that we obtain the consistency with the discrete continuity equation in the sense of \cite{gross:2002} by using the values of the height at time level $n$, $h^n$, in the advection terms. In the equations above we use the discrete transference term
$$
\begin{array}{lll}
G^{n,l}_{\a+\frac12,i} &=&  \dfrac{1}{\dxi}\dsum_{\b=1}^{\a} \Bigg(  l_{\b,i+\frac12}h^n_{i+\frac12}u^{n,l}_{\b,i+\frac12} \,-\,l_{\b,i-\frac12}h^n_{i-\frac12}u^{n,l}_{\b,i-\frac12}    \\
&-&  l_{\b,i}\dsum_{\gamma=1}^{N}\left(  l_{\gamma,i+\frac12}h^n_{i+\frac12}u^{n,l}_{\gamma,i+\frac12} \,-\,l_{\gamma,i-\frac12}h^n_{i-\frac12}u^{n,l}_{\gamma,i-\frac12}    \right)\Bigg),
\end{array}
$$ and $\r^{n,l}_{\a+\frac12,i}$ is the upwind value depending on the vertical velocity. To this aim, we define
\begin{equation}
\label{eq:disc_dens_equation}
\r_{\a+1/2,i}^{n,l}=\frac{\r^{n,l}_{\a,i}+\r^{n,l}_{\a+1,i}}{2} \ -\  \dfrac{\mbox{sgn}\left(-G^{n,l}_{\a+\frac12,i}\right)}{2}\left(\r^{n,l}_{\a+1,i}-\r^{n,l}_{\a,i}\right),
\end{equation}
for $l=1,2,3$.

 \section{Numerical results}
\label{se:tests}

We present in this section some numerical tests in order to validate the proposed method through some academic configurations of variable density flows. Based on the linear analysis in section \ref{se:linear}, we define a Courant number
based on celerity as

\begin{equation}\label{eq:cfl_cel}
	C_{cel} = \displaystyle\max_{1 \leq i \leq M}\, \left(\abs{\bar{u_{i}}} + \sqrt{\left(1+ \bar \r_{i}\right)g\,h_i}\right)\,\dfrac{\dt}{\dxi},\quad
	 \bar {u} = \dsum_{\a=1}^N l_\a u_\a,\quad 
	\bar{\r} = \dsum_{\a=1}^N l_\a \r_\a,
	\end{equation}
	where  $\bar{u}, \bar{\r}$ are the vertically averaged velocity and density perturbation corresponding to each control volume, and a Courant number based on velocity  as
\begin{equation}\label{eq:cfl_vel}
	C_{vel} = \displaystyle\max_{1 \leq i \leq M}\, \left(\abs{\bar{u_{i}}} + \sqrt{\bar \r_{i}g\,h_i}\right)\,\dfrac{\dt}{\dxi},\quad
	 \bar {u} = \dsum_{\a=1}^N l_\a u_\a,\quad 
	\bar{\r} = \dsum_{\a=1}^N l_\a \r_\a,
	\end{equation}
Since the terms associated to the barotropic pressure gradient are treated implicitly in our approach, the resulting stability
condition will be  based on $C_{vel} $ rather than $C_{cel}, $ thus allowing for substantial computational gains.

For all the tests, we   compute a reference solution   with a third-order Runge-Kutta method with a fixed value $C_{cel}=0.1$. That is, the time-step is adaptive, as a function of $C_{cel}$ \eqref{eq:cfl_cel}.  
The accuracy is measured with the following definitions: $Err_\eta\,[\,l_2\,]$ and  $Err_\eta\,[l_\infty\,]$ denote the relative error for the free surface, and for the velocity and density fields we use
\begin{subequations}\label{eq:error}
	\begin{alignat}{2}
	Err_s[l_2] &= \left(\dfrac{\sum_{\a = 1}^N \sum_{i = 1}^M \abs{s_{\a,{i+\frac12}}-s_{\a,{i+\frac12}}^{ref}}^2\dxi h_{\a,i}  }{\sum_{\a = 1}^N \sum_{i = 1}^M \abs{s_{\a,{i+\frac12}}^{ref}}^2\dxi h_{\a,i} }\right)^{1/2};\\[4mm]
	Err_s[l_\infty] &= \dfrac{\max_{\a}\max_{i} \abs{s_{\a,{i+\frac12}}-s_{\a,{i+\frac12}}^{ref}} }{\max_{\a}\max_{i} \abs{s_{\a,{i+\frac12}}^{ref}}  },
	\end{alignat}
\end{subequations}
where $s=\r$ or $s=u$, and $s^{ref}$ denotes the reference solution.

Note that a factor depending on the density is included in the Courant numbers \eqref{eq:cfl_cel} and \eqref{eq:cfl_vel}. In these definitions, 
$\bar{\r}$ is used to obtain an approximation of the gravity wave speed associated to the density perturbation. 
 Finally, all the computational times showed in this section have been measured on a Mac Mini with Intel$\textsuperscript{\textregistered}$Core$\texttrademark$ i7-4578U and 16 GB of RAM.

\subsection{Internal gravity wave}
\begin{figure}[!th]
	\centering\includegraphics[width=0.85\textwidth]{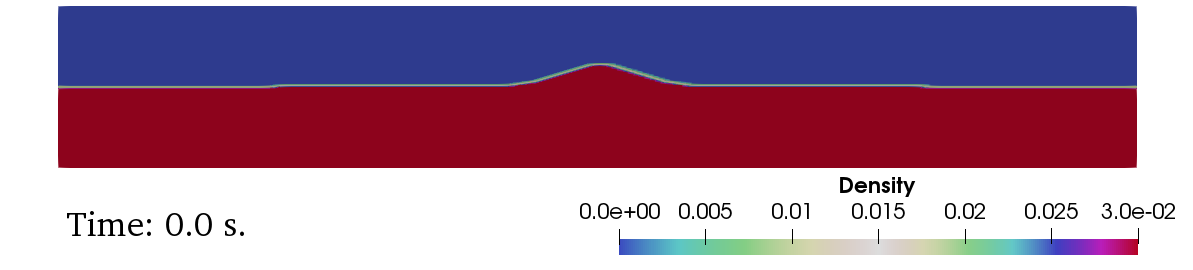}
	\centering\includegraphics[width=0.85\textwidth]{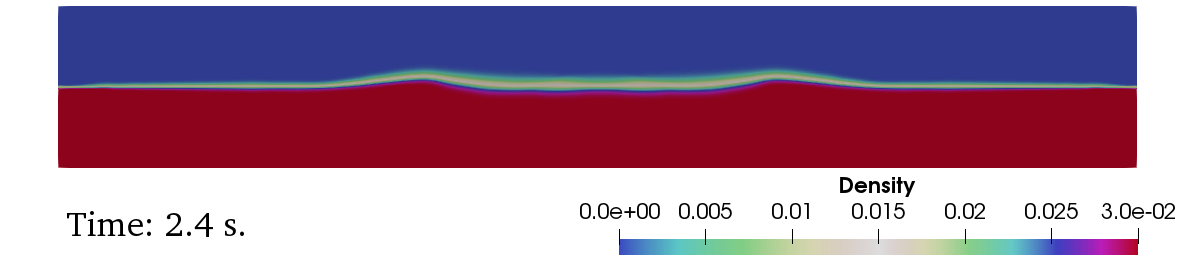}
	\centering\includegraphics[width=0.85\textwidth]{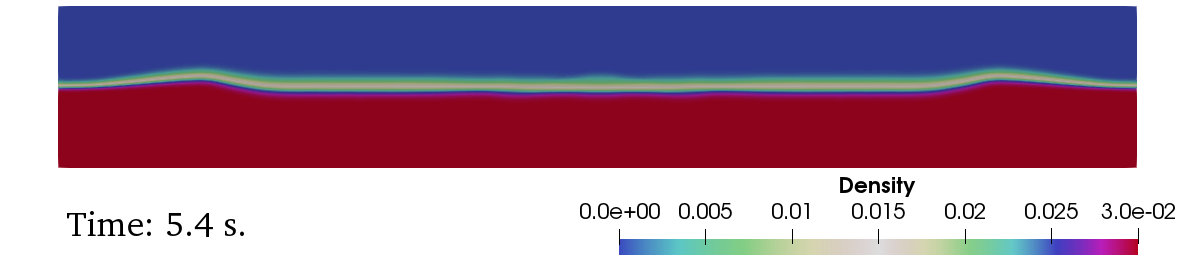}
	\caption{\label{fig:onda_interna} \footnotesize \textit{Density field for different times: $t = 0,\,2.4$ and $5.4$ s.}}
\end{figure}

We consider here a internal gravity wave produced by a perturbation in the density field. The computational domain is $[0,L]$, with $L=2$ m. It is supposed to be closed (wall boundary conditions), with flat bottom, and the fluid at rest at initial time, when the initial height is $\eta_0(x) = 0.3$ m everywhere. The fluid have two separate areas with densities $\r_1 = 1000$ kg/m$^3$ and $\r_2 = 1030$ kg/m$^3$. These values are rewritten in terms of the relative perturbation as $\r_1 = 0$ and $\r_2 = 0.03$. Thus, the initial condition for the density is given, for $\a=1,\dots,N$, by   
$$
\r_{0,\a}(x) = \left\{ \begin{array}{ll}
0.03 & \mbox{if }z_\a < z_{lim}\\
0 & \mbox{otherwise} 
\end{array}\right.$$
where $z_{lim} = 0.15 + 0.04 e^{-100\left(x-1\right)^2}$ and $z_\a = \eta_0 \dsum_{i=1}^\beta l_\beta$ (see figure \ref{fig:onda_interna}).
We take $\dx=0.01$ m and a non-uniform distribution of the vertical layers in order to get an accurate definition of the density perturbation. To this aim, we consider that the layers are concentrated over the central part (along the vertical direction) of the domain. Then,  $54$ layers are used, where four of them are in the first and last quarter of the domain, and $50$ layers over the central part. Therefore, the coefficients $l_\a$ are defined as
$$
l_1 = l_2 = 0.125,\quad l_i = 0.01, i=3,\dots,52,\quad l_{53}=l_{54} = 0.125.
$$
Note that this definition of the layers should be equivalent to consider $100$ uniform layers. Actually, we have checked that obtained results are similar in both cases.

\begin{figure}[!th]
	\centering\includegraphics[width=0.55\textwidth]{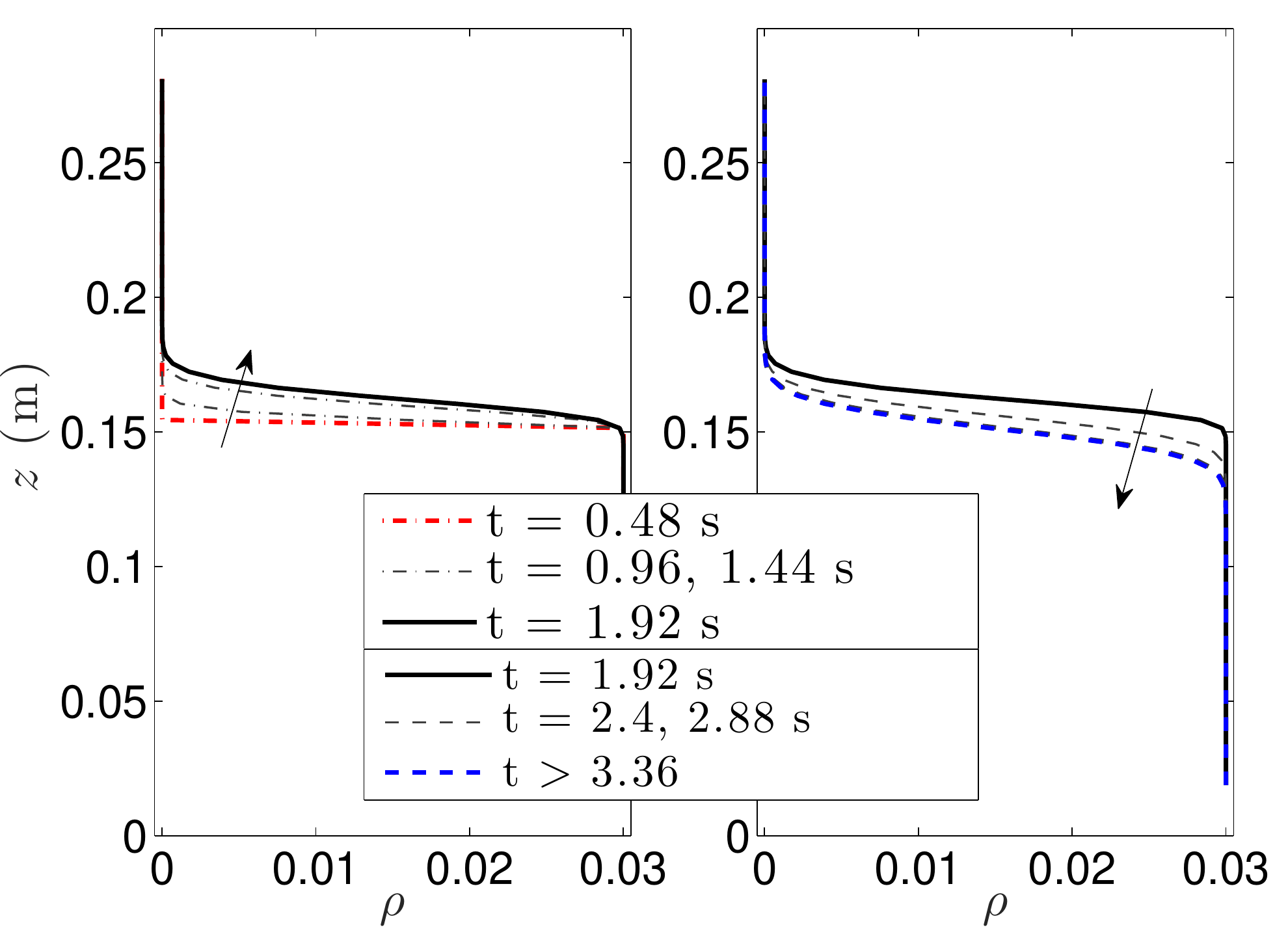}
	\caption{\label{fig:onda_perfil} \footnotesize \textit{Time evolution of vertical profiles of density at the point $x=0.745$ m.}}
\end{figure}

In figure \ref{fig:onda_interna} we see the density field for the initial time, and also for $t=2.4,\ 5.4$ s. As expected, the initial density perturbation is subdivided in two internal waves, which travel in opposite directions. In this test, there is no significant variation of the initial height. Figure \ref{fig:onda_perfil} shows the time evolution of the density at the point $x=0.745$ m. We see one of the generated waves getting through that point. Thus, the interface between the areas with different density rises and then falls. In these figures we see the numerical dissipation introduced by the vertical discretization. It could be improved by using a more accurate discretization of the equation for the density perturbation, namely for the definition of $\r_{\a+1/2}$ given by \eqref{eq:disc_dens_equation}. Nevertheless, it is not the goal of this paper but showing how the proposed semi-implicit method is efficiently adapted to variable density flows. 

Table \ref{tab:tabla_error_internalWave} shows the relative errors and Courant numbers for this test, when comparing with the solution computed using the explicit method, at time $t=4.8$ s, before the perturbations arrive to the boundaries. The errors for the velocity field are greater than the ones for the density field by one order of magnitude. We have reasonable errors till Courant number $C_{cel} = 7.0$, which are lower than $8\%$ and $1\%$ for the velocity and density fields, respectively. For $C_{cel} = 13.8$, these errors are $10\%$ approximately for the velocity. When a larger time step is used, the errors grow up quickly and the simulation becomes unstable for $\dt \approx 0.09$, corresponding to $C_{vel} \approx 2.04$. 


\begin{table}[!h]
	\begin{center}
		\begin{tabular}{cccccccc}
			\hline
			$\dt$ (s)   & $C_{cel}$ & $C_{vel}$ & Err$_{\eta} $ [$l_2/l_\infty$]  & Err$_{u}$ [$l_2/l_\infty$]& Err$_{\rho}$ [$l_2/l_\infty$] \\ 
			& & &    ${\small(\times10^{-4})}$ &  ${\small(\times10^{-2})}$&  ${\small(\times10^{-2})}$\\ \hline\\
			0.01 &  1.7 & 0.22 & 0.8/2.3  & 2.9/1.5 & 0.03/0.1\\
			0.02 &  3.5 & 0.45 & 0.7/1.8  & 7.7/6.9 & 0.2/1.1\\
			0.04 &  7.0 & 0.91 & 2.3/5.4  & 7.3/15.2 & 0.9/6.4\\
			0.06 &  10.3 & 1.35 &  1.8/3.9 & 10.4/27.4 & 1.5/10.8\\
			0.08 & 13.8  & 1.81 & 1.8/3.7  & 10.5/21.7 & 1.8/14.3\\
			\hline\end{tabular}
		\caption{\footnotesize \it{ Relative errors and Courant numbers achieved by using the IMEX-ARK2 method at $t= 4.8$ s. }}
		\label{tab:tabla_error_internalWave}
	\end{center}
\end{table}

Table \ref{tab:tabla_speedup_internalWave} shows the computational times and speed-up that we obtain with the semi-implicit methods for a final time $t_f = 10$ s. We see that with $\dt = 0.04$ s, the semi-implicit method is almost $8$ times faster than the explicit method with $C_{cel} = 0.9$, and for $\dt = 0.08$ the speed-up is $16$.

\begin{table}[!h]
	\begin{center}
		\begin{tabular}{ccccccccc}
			\hline\\
			Method & $\dt$ (s)  & $C_{cel}$ & Comput. time (s)  & Speed-up\\ 
			& & &    &  \\ \hline\\
			Runge-Kutta 3 & - & 0.1 & 719.3 (12 min.) & -\\\
			Runge-Kutta 3 & - & 0.9 & 83 & 1\\
			IMEX-ARK2 & 0.01 & 1.7 & 42.1 & 2.0\\
			IMEX-ARK2 & 0.02 & 3.5  & 20.8 & 4.0\\
			IMEX-ARK2 & 0.04 & 7.0 &  10.8 & 7.7\\
			IMEX-ARK2 & 0.06 &  10.3 & 7.1 & 11.7\\
			IMEX-ARK2 & 0.08 &  13.8 & 5.2 & 15.9\\
			\hline\end{tabular}
		\caption{\footnotesize \it{ Speed-ups achieved with the IMEX-ARK2 method at final time $t_f= 10$ s.}}
		\label{tab:tabla_speedup_internalWave}
	\end{center}
\end{table}

%

\subsection{Lock exchange}
\begin{figure}[!h]
	\centering\includegraphics[width=0.85\textwidth]{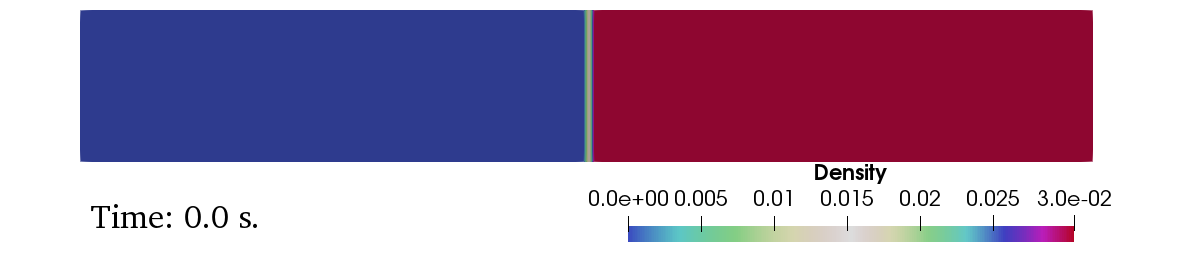}
	\centering\includegraphics[width=0.85\textwidth]{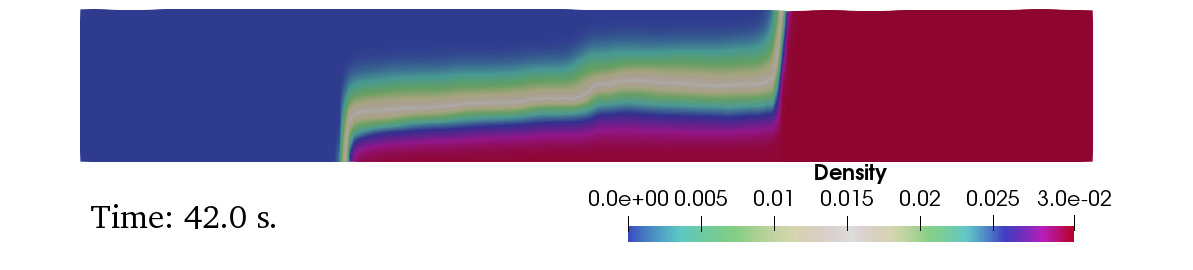}
	\centering\includegraphics[width=0.85\textwidth]{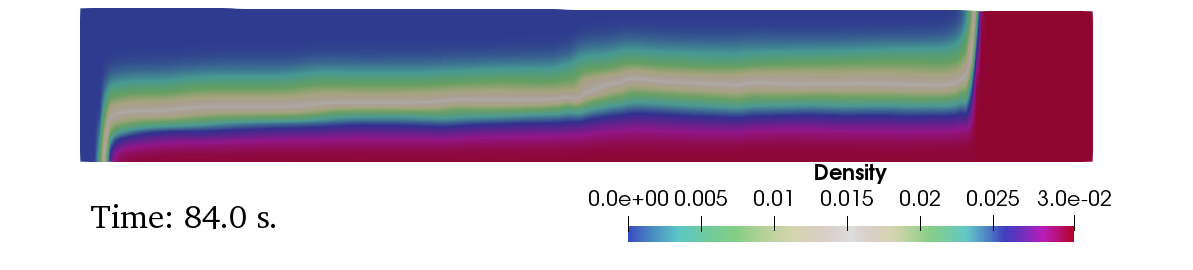}
	\caption{\label{fig:lock_exchange} \footnotesize \textit{Density distribution for initial time and times $t= 42,84$ s.}}
\end{figure}

Let us consider now a classical test for variable density flows, the lock exchange problem, where the fluid that is located on the right hand side of the domain have a density higher than the fluid located on the left hand side. We consider a closed domain, with $x\in[-L/2,L/2]$, whose length is $L=20$ m and $\dx=0.1$ m. In this case we take $20$ uniform layers to reproduce the vertical structure of the density field.

The flow starts from the rest, with initial height $\eta_0(x) = 0.3$ m, while the initial density is defined by the function
$$
\r_{0,\a}(x) = \left\{ \begin{array}{ll}
0.03 & \mbox{if } x>0,\\
0 & \mbox{otherwise}, 
\end{array}\right.\qquad \mbox{for }\a=1,\dots,N.$$

In this case, the discontinuous initial profile of density makes necessary to apply a flux limiter to the second-order discretization of the advection term $u\p_x u$. Otherwise, spurious oscillations appear in the simulation. Here we consider  the classical {\it minmod} flux limiter (see e.g. \cite{leveque:2002}).

The behaviour is the expected one, qualitatively. Figure \ref{fig:lock_exchange} shows the density field at initial time, an intermediate time ($t=42$ s), and the time where we measure the errors between explicit and the semi-implicit method ($t=84$ s). 

Table \ref{tab:tabla_error_lockExchange} shows the relative errors and Courant numbers achieved, where we got a Courant number $C_{cel} = 5.3$. Even thought the results are good enough, we see that these errors are greater than in previous test, namely for the density field. In this case, the error for the velocity and density fields are lower than $3\%$ for $C_{cel} = 3.5$. However, high values of the errors are observed, in particular using the norm $l_{\infty}$, when $C_{cel}$ are larger. It is due to the fact that the vertical structure of the fluid is stronger in this test,
therefore the discretization of the mass transference terms start playing a role and 
only allows to achieve Courant numbers lower than for the previous test. 
In addition, as commented before, since the goal of this work is to evaluate the accuracy of the proposed semi-implicit method, we have chosen a first order discretization of the equation of the density evolution \eqref{eq:dens_var}. If one need to reduce the error for the density and velocity, a more accurate discretization of that equation could be used. 

\begin{table}[!h]
	\begin{center}
		\begin{tabular}{ccccccccc}
			\hline\\
			$\dt$ (s) &  $C_{cel}$ & $C_{vel}$ & Err$_{\eta} $ [$l_2/l_\infty$]  & Err$_{u}$ [$l_2/l_\infty$]& Err$_{\rho}$ [$l_2/l_\infty$] \\ 
			& &  &   ${\small(\times10^{-3})}$ &  ${\small(\times10^{-2})}$&  ${\small(\times10^{-2})}$\\ \hline\\
			0.1  & 1.8 & 0.31 & 0.6/1.8  & 1.2/2.2 & 0.3/2.5\\
			0.2  & 3.5 & 0.62 & 0.8/1.9  & 2.7/15.0 & 1.3/13.0\\
			0.3 & 5.3 & 0.93 & 1.3/3.2  & 7.5/50.0 & 4.8/56.8\\
			\hline\end{tabular}
		\caption{\footnotesize \it{ Relative errors and Courant numbers achieved by using the IMEX-ARK2 method at $t= 84$ s.}}
		\label{tab:tabla_error_lockExchange}
	\end{center}
\end{table}

In table \ref{tab:tabla_speedup_lockExchange} we see the speed-ups achieved for the semi-implicit method. We get a speed-up of 5.3 in this case.

\begin{table}[!h]
	\begin{center}
		\begin{tabular}{ccccccccc}
			\hline
			Method & $\dt$ (s)  & $C_{cel}$ & Comput. time (s)  & Speed-up\\ 
			& & &    &  \\ \hline
			Runge-Kutta 3 & - & 0.1 & 197.11 (3.3 min.) & -\\
			Runge-Kutta 3 & - & 0.9 & 21.3 & 1\\
			IMEX-ARK2 & 0.1 & 1.8 & 12.0 & 1.8\\
			IMEX-ARK2 & 0.2 & 3.5 & 6.0 & 3.6\\
			IMEX-ARK2 & 0.3 & 5.3 & 4.0 & 5.3\\
			\hline\end{tabular}
		\caption{\footnotesize \it{ Speed-ups achieved with the IMEX-ARK2 method at final time $t_f= 100$ s.}}
		\label{tab:tabla_speedup_lockExchange}
	\end{center}
\end{table}

We have also measured the front velocity in both, the surface and the bottom. Following  \cite{hervouet:2007}, one could estimate that, if all the potential energy in the initial condition is transformed into kinetic energy, the theoretical velocity of the front should be
$V = \sqrt{0.25 g h \r},$
where $\r$ is the density perturbation. In this case, it leads to $V = 0.1485$ m/s. 
We have also numerically computed the front velocity in both the surface and the bottom. We obtain identical values for the explicit and semi-implicit methods at any time step. These velocities, measured at time $t=10$ s, are $0.145$ m/s for the surface and $0.095$ m/s for the bottom, which leads to a mean error of $19\%$, approximately.

\subsection{Tidal forcing with variable density}\label{test:marea}
Now, we look for a more realistic flow, simulating the mouth of a river into the sea. We consider the computational domain $x\in [-7500,22500]$ m, whose length is $30$ km, and the variable bathymetry is given by 
$$b(x) = z_0 - z_1\,\text{tanh}(\lambda\,(x-x_0)) +  20\,e^{-(x-x_1)^2/\sigma^2},$$
with $z_0  = -z_1 = 46$, $x_0 = 7500$, $x_1 = 16000$, $\lambda = -1/3000$ and $\sigma=2000$.  This test is analogous to the one in \cite{bonaventura:2018b} with some differences. Mainly, the height of the bump in the bottom is $20$ m. We take a smaller pick in order to ensure the subcritical regime, since an hydraulic jump due to the variable density is observed with the bottom definition in \cite{bonaventura:2018b}. Also the height of the shallowest part of the domain is lower, in order to properly reproduce the saltwater intrusion into a river. 

We consider a fluid with constant density, $\rho_0 = 1000$ kg/m$^3$, in the whole domain at initial time. Downstream the water going into the domain have a higher density $\rho_1 = 1030$ kg/m$^3$. To this aim, we assume that the water is at rest at initial time, and the deviation of the constant density is zero everywhere ($\r_{0,\a}(x) = 0$). The initial free surface is defined by $\eta_0(x) = 100$ m. We consider $10$ uniform vertical layers, and $500$ nodes in the horizontal discretization, i.e. $\dx = 60$ m. Now, subcritical boundary condition are imposed:
\begin{itemize}
	\item Downstream: a tidal downstream condition $\eta(t,L) = 100 + 3\sin(\omega t)$ m is assumed, where $\omega = 2\pi/43200$. For the density, in order to avoid spurious oscillation appearing as consequence of the discontinuous boundary condition, we account for the vertical structure of the flow. In particular, we run a simulation with constant density perturbation $\r_{\a}(t,L) = 0.03$ for $\a = 1,\dots,N$, and after some periods of tide (when the dynamics of the flow is stabilized), the profile of the density perturbation close to the boundary is $\r_{\a,ext} = 0.03, \a = 1,...,7$, and
	$$
	\r_{8,ext} = 0.028,\ \r_{9,ext} = 0.025,\ \r_{10,ext} = 0.015.
	$$
	Then, we consider as boundary condition the obtained profile $\r_{\a}(t,L) = \min \left(\r_{\a,ext},\r_{\a,ext} t/(6\cdot 3600)\right)$, that is, the density perturbation goes into the domain slowly and smoothly.
	
	\item Upstream: for the discharge, we define $q(t,-7500) = \min \left(1,t/(6\cdot 3600)\right)$ m$^2$/s, and we impose fresh water going into the domain $\r_\a(t,-7500) = 0$, for $\a = 1,\dots,N$. 
\end{itemize}

In this case we also consider a turbulent viscosity as in \cite{bonaventura:2018b} (we refer to previous work for details), with friction parameters $\D z_r = l_1 h$ m, $\D z_0 = 3.3 \times10^{-5}$ m and $\kappa = 0.41$. The wind drag is defined by $C_w = 1.2\times 10^{-6}$ and wind velocity $u_{w}=1$ m/s. Notice that, for the sake of simplicity, we consider this turbulence model which does not take into account the density perturbation. In realistic applications, more accurate turbulence models would have to be applied.

\begin{figure}[!ht]
	\centering\includegraphics[width=1\textwidth]{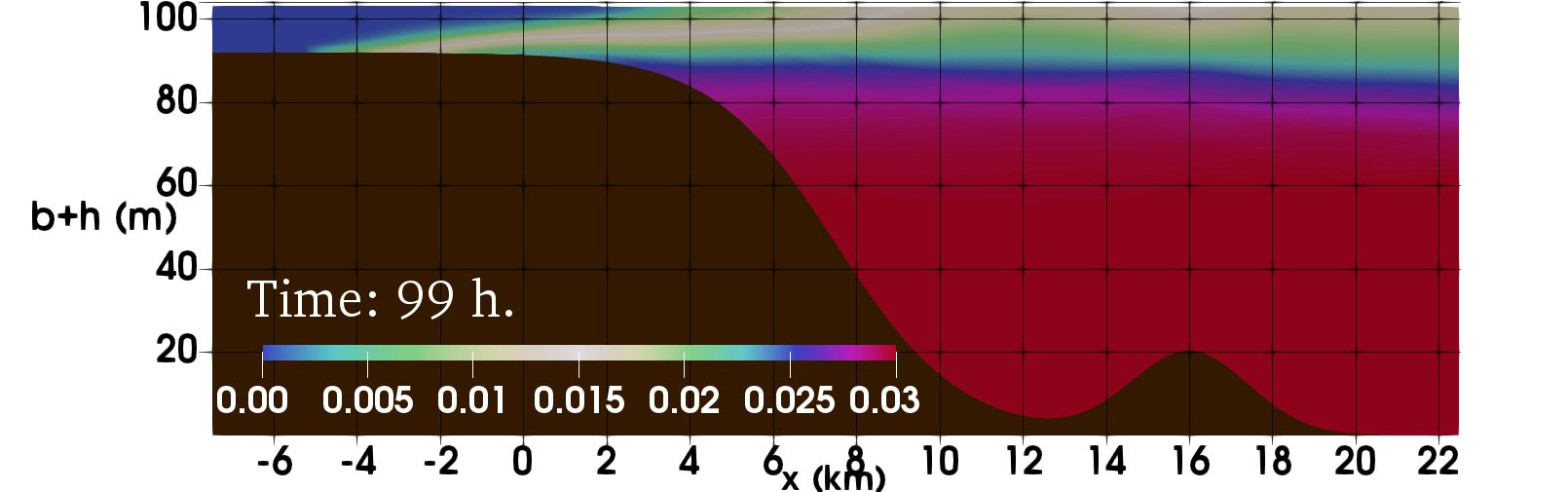}\\[3mm]
	\centering\includegraphics[width=1\textwidth]{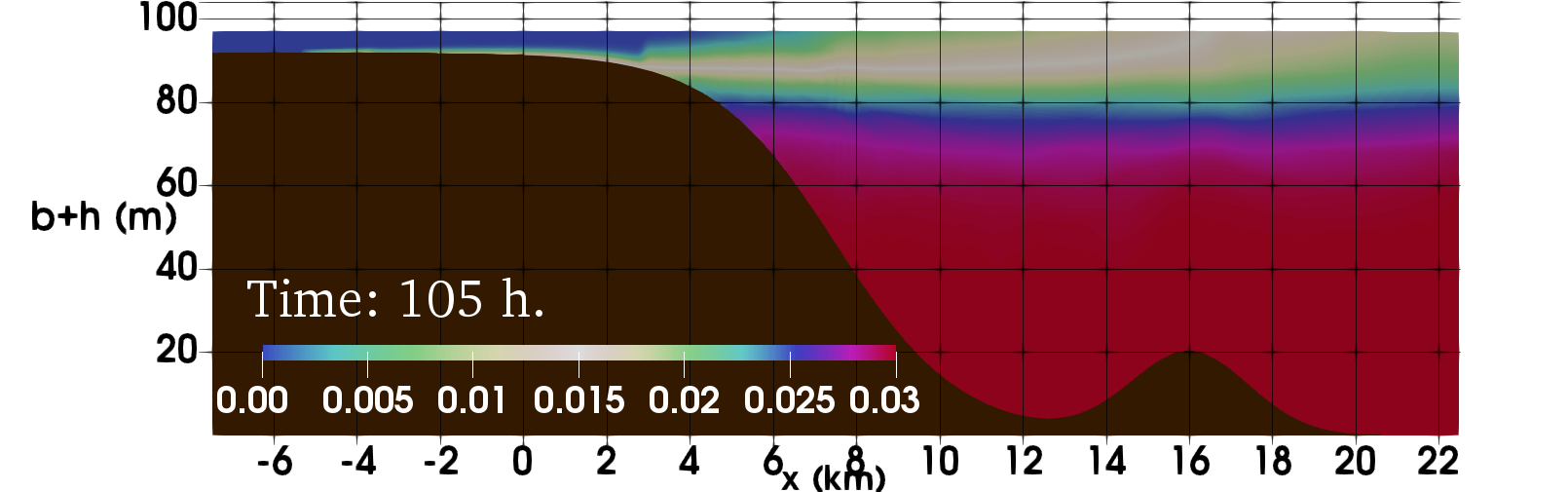}\\[3mm]
	\centering\includegraphics[width=1\textwidth]{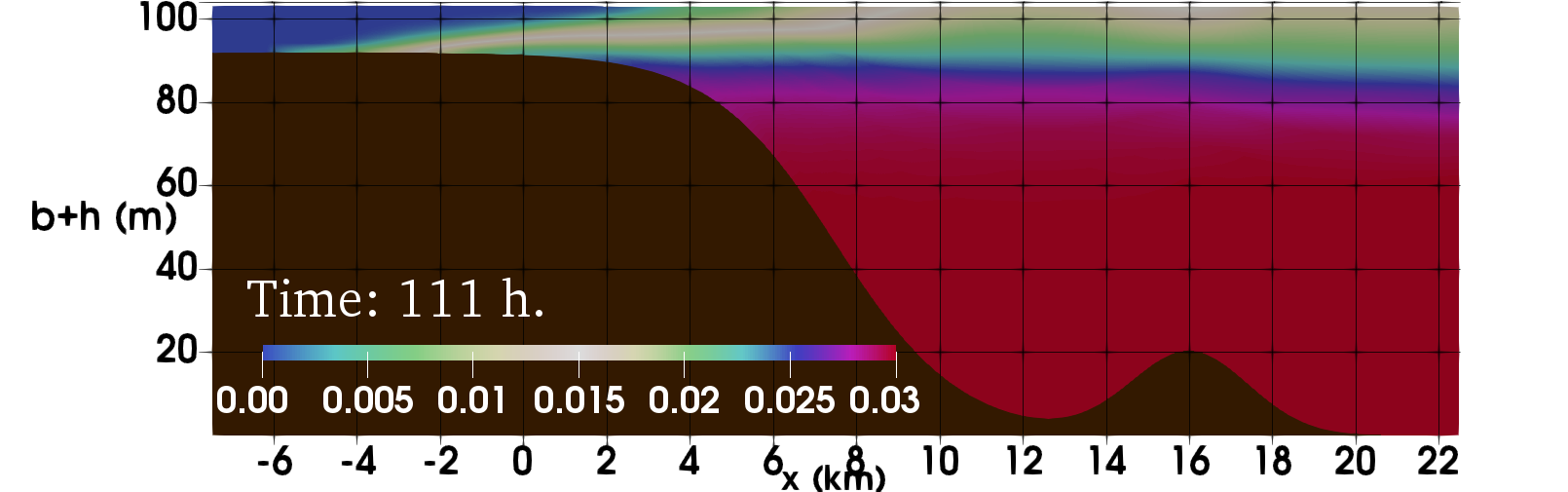}\\[3mm]
	\centering\includegraphics[width=1\textwidth]{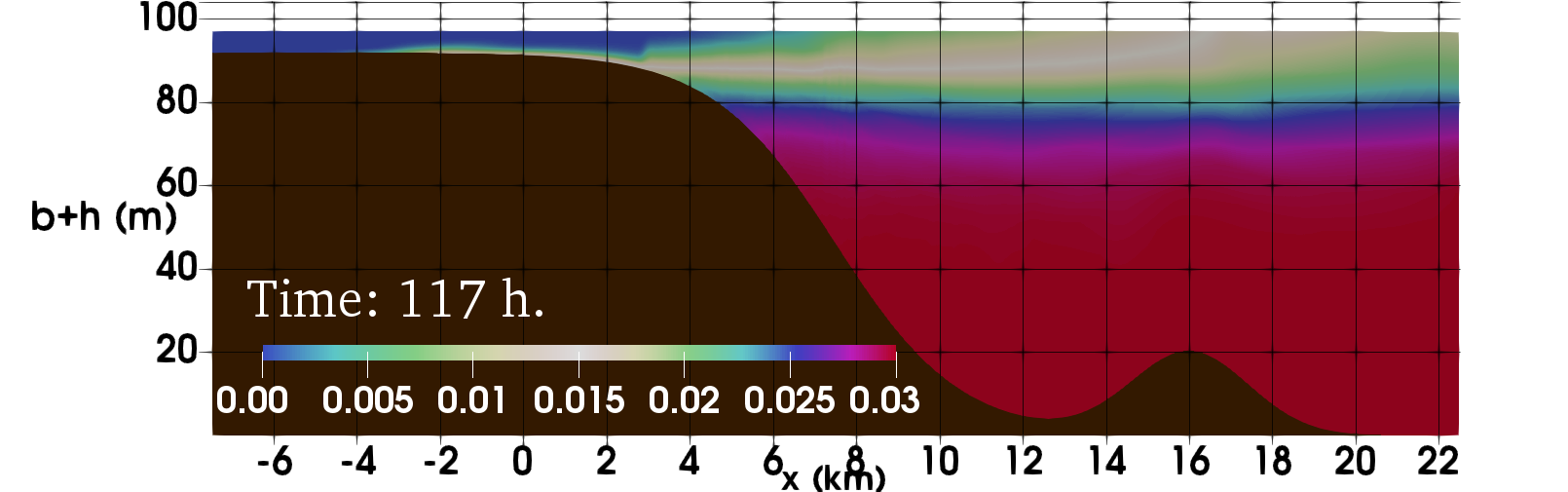}
	\caption{\label{fig:marea_density2} \footnotesize \textit{Saltwater intrusion into the river at times $t=99,105,111,117$ hours.}}
\end{figure}

\begin{figure}[!th]
	\begin{minipage}{16cm}
		\hspace{-3cm}	\centering\includegraphics[width=0.49\textwidth]{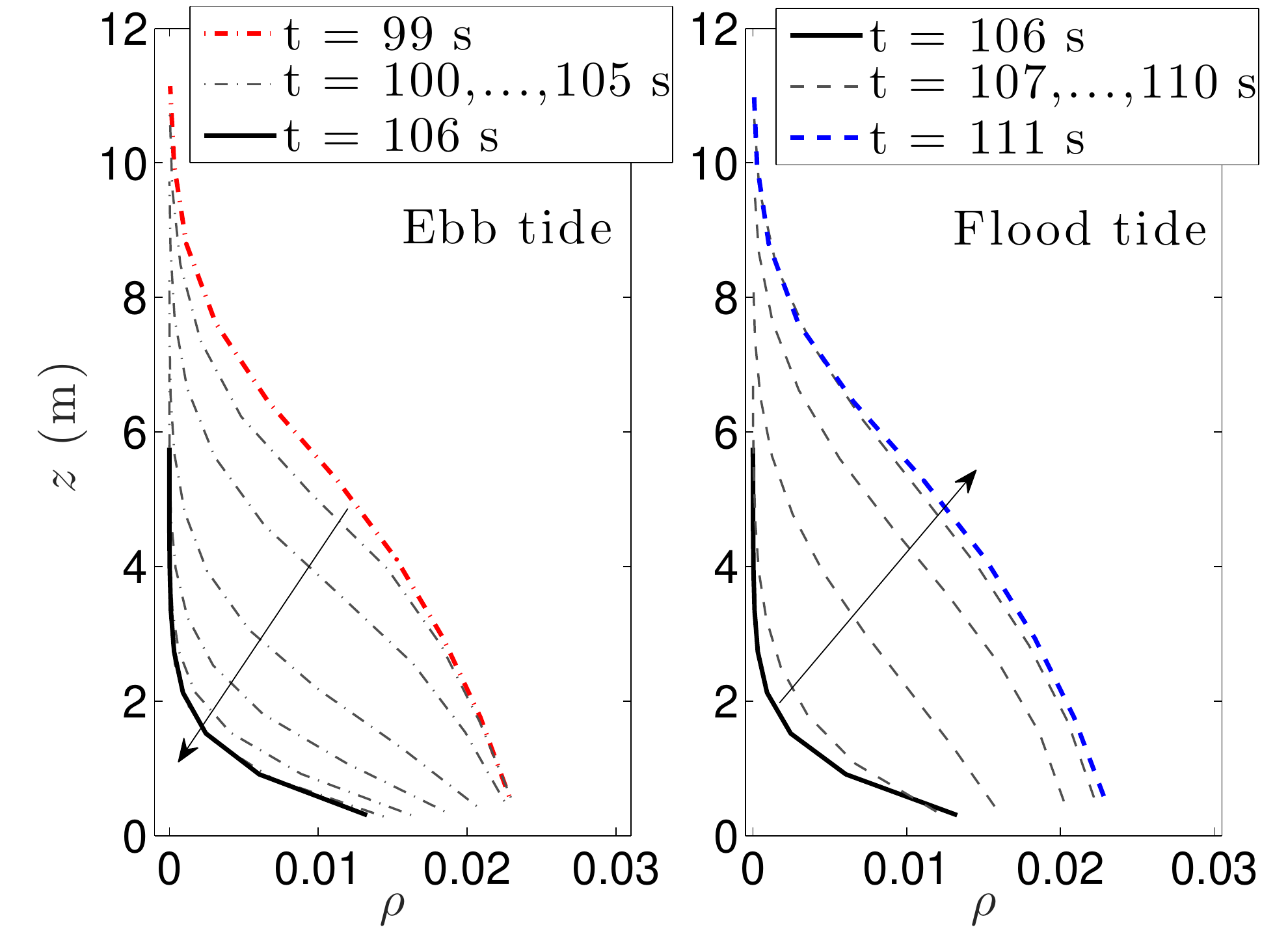}
		\centering\includegraphics[width=0.49\textwidth]{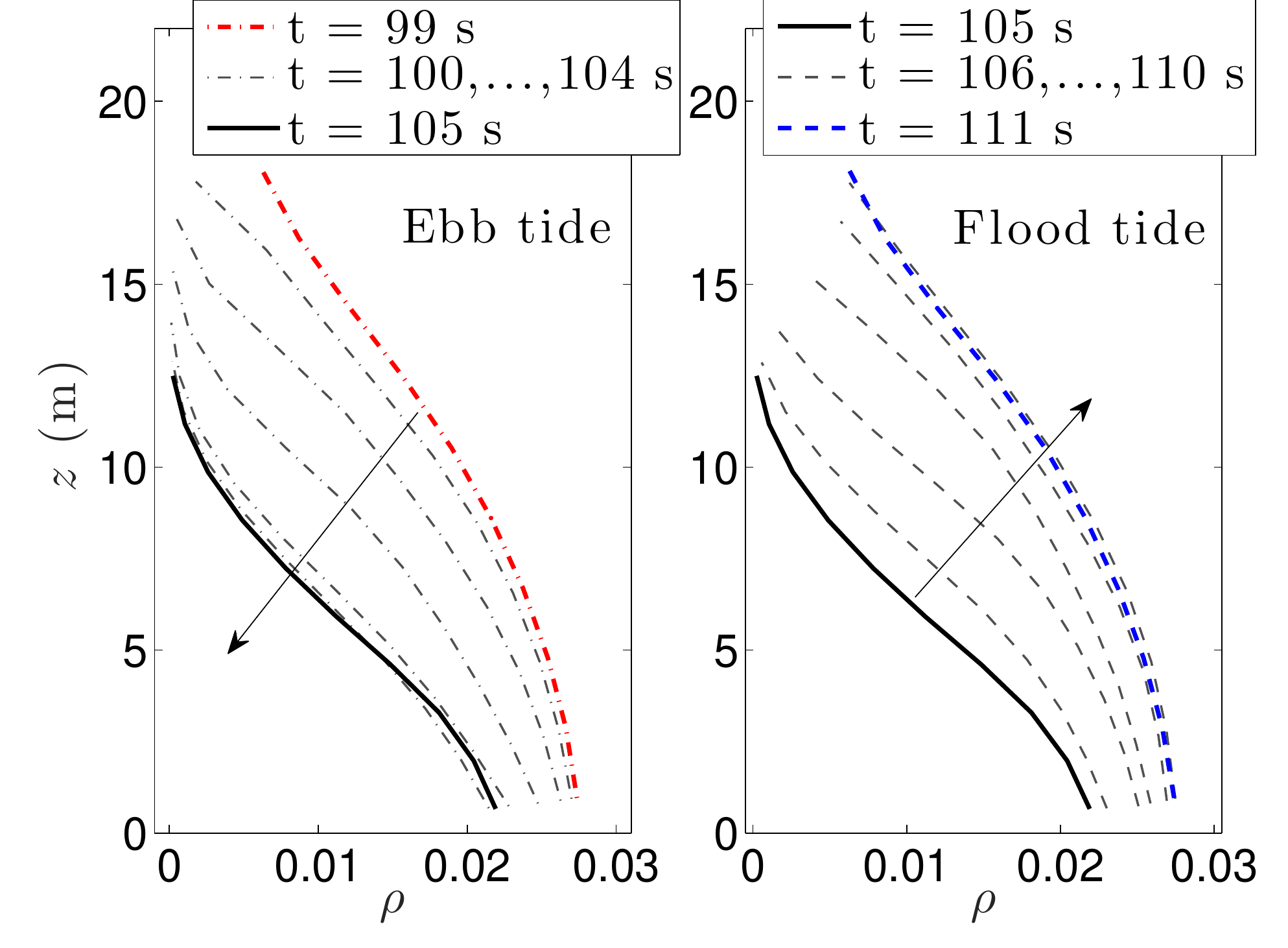}
	\end{minipage}
	\caption{\label{fig:marea_perfil} \footnotesize \textit{Time evolution of vertical profiles of density at the point $x=0$ km (left hand side) and $x=4$ km (right hand side)}}
\end{figure}

We simulate twelve 12-hours periods of tide, i.e., 144 h. Figure \ref{fig:marea_density2} shows an example of $1.5$ periods of tide. The density distribution is the expected one, the water with higher density goes to the low part of the domain, and a layer with lower density is observed in an upper layer, especially in the right-hand side (deepest area) of the domain. The saltwater intrusion phenomenon is well reproduced, and we obtain the expected periodic behaviour. It is also observed in Figure \ref{fig:marea_perfil}, where the vertical profiles of density are shown at points $x = 0$ km and $x=4$ km (mouth of the river). We see the water level falling (ebb tidal phase), and most of the water column is fresh water. Next, the water level rises (flood tidal phase) and saltwater is going into the river.

\begin{table}[!h]
	\begin{center}
		\begin{tabular}{ccccccccc}
			\hline
			$\dt$ (s) &  $C_{cel}$ & $C_{vel}$ & Err$_{\eta} $ [$l_2/l_\infty$]  & Err$_{u}$ [$l_2/l_\infty$]& Err$_{\rho}$ [$l_2/l_\infty$] \\ 
			& &  &   ${\small(\times10^{-6})}$ &  ${\small(\times10^{-2})}$&  ${\small(\times10^{-2})}$\\ \hline
			5  & 2.7 & 0.45 & 4.1/12.5  & 1.0/4.8 & 0.3/2.1\\
			10  & 5.4 & 0.9 & 3.0/8.5  & 0.8/4.1 & 0.2/0.8\\
			15 & 8.1 & 1.35 & 3.7/11.2  & 1.0/2.9 & 0.2/1.2\\
			20 & 10.7 & 1.8 & 4.1/14.9  & 1.2/7.4 & 0.2/1.4\\
			\hline\end{tabular}
		\caption{\footnotesize \it{Relative errors and Courant numbers achieved by using the IMEX-ARK2 method at $t= 144$ h.}}
		\label{tab:tabla_error_marea}
	\end{center}
\end{table}

In Table \ref{tab:tabla_error_marea} we see the relative errors and Courant numbers achieved at final time $t=144$ h. We have a Courant number $C_{cel} = 10.7$ with $C_{vel} = 1.8$. The $l_2$-errors in the free surface position has order $10^{-6}$, while the errors are larger for the velocity and density fields. 

\begin{table}[!h]
	\begin{center}
		\begin{tabular}{ccccccccc}
			\hline
			Method & $\dt$ (s)  & $C_{cel}$ & Comput. time (min)  & Speed-up\\ 
			& & &    &  \\ \hline
			Runge-Kutta 3 & - & 0.1 & 597.3 (9.95 h) & -\\\
			Runge-Kutta 3 & - & 0.9 & 68.31 & 1\\
			IMEX-ARK2 & 5 & 2.7 & 28.95 & 2.4\\
			IMEX-ARK2 & 10 & 5.4 & 14.52 & 4.7\\
			IMEX-ARK2 & 15 & 8.1 & 9.69 & 7.0\\
			IMEX-ARK2 & 20 & 10.7 & 7.25 & 9.4\\
			\hline\end{tabular}
		\caption{\footnotesize \it{Speed-ups achieved with the IMEX-ARK2 method at final time $t_f= 144$ h.}}
		\label{tab:tabla_speedup_marea}
	\end{center}
\end{table}

Table \ref{tab:tabla_speedup_marea} shows the speeds-up achieved with the semi-implicit method for a final time $t=144$ hours.  The IMEX discretization is almost $5$ times faster for $\Delta t = 10$ s ($C_{vel}=0.9$), and $11$ times faster for $\Delta t = 20$ s.

In the following, we analyze the possibility os reducing the number of degrees of freedom of the multilayer system in the shallowest part of the domain. We study several configurations trying to reduce the error made, and preserving the vertical structure if needed. They are denoted hereinafter as (NVAR(n)), where (n) indicates the number of layers used in each configuration. Thus, we simplify the vertical discretization in the shallowest part of the domain ($x<0$ km) as follows:
The vertical discretization is totally removed and a single layer is considered in the first part of the domain:
\begin{equation}\label{eq:marea_nvar1}
\tag{NVAR(1)}
N = \left\{\begin{array}{lll}	
1, & l_1 = 1, & \mbox{ if } x\leq 0;\\
10, & l_i = 1/10, \,i=1,...,N,  & \mbox{ otherwise}.\\
\end{array}\right.
\end{equation}
The following configurations correspond to non-uniform distribution of the layers. Three and four layers are used in the shallowest part of the domain,  increasing also the thickness of the layers close to the bottom in order to improve the vertical discretization, namely the friction effect:
\begin{equation}\label{eq:marea_nvar3}
\tag{NVAR(3)}
N = \left\{\begin{array}{lll}	
3, & l_1 = l_2 = 0.2, l_3 = 0.6, & \mbox{ if } x \leq 0;\\
10, & l_i = 1/10, \,i=1,...,N,  & \mbox{ otherwise},\\
\end{array}\right.
\end{equation}
and
\begin{equation}\label{eq:marea_nvar4}
\tag{NVAR(4)}
N = \left\{\begin{array}{lll}	
4, & l_1 = l_2 = 0.2, l_3 = 0.2,l_4=0.4, & \mbox{ if } x \leq 0;\\
10, & l_i = 1/10, \,i=1,...,N,  & \mbox{ otherwise}.\\
\end{array}\right.
\end{equation}
For these configurations with a variable number of vertical layers, we consider the IMEX method with $\Delta t = 10$ s in all the simulation.

\begin{figure}[!ht]
	\begin{minipage}{16cm}
		\hspace{-3cm}	\centering\includegraphics[width=1.\textwidth]{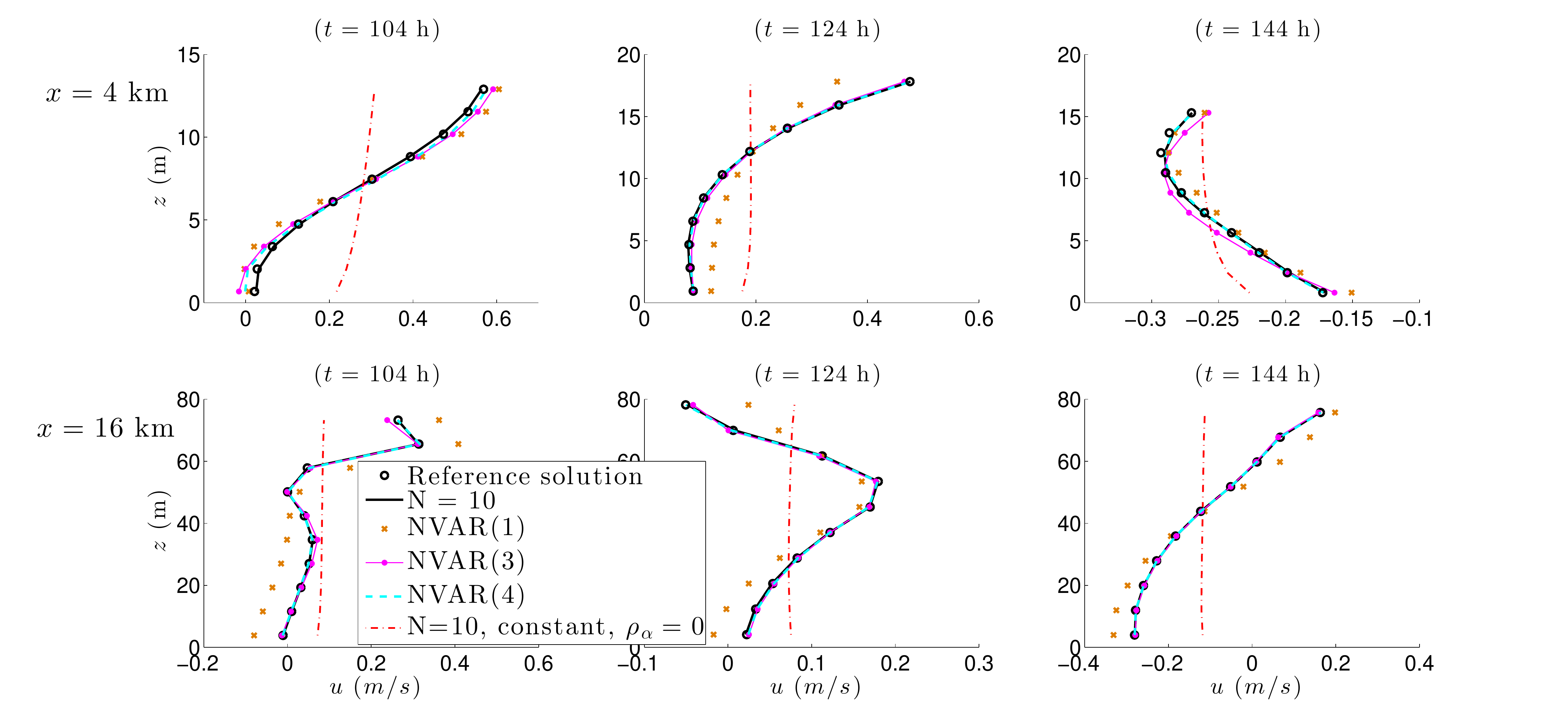}
	\end{minipage}
	\caption{\label{fig:marea_nvar_vel} \footnotesize \textit{Vertical profiles of horizontal velocity obtained with the reference solution (black circles), and the IMEX with $\dt = 10$ s, and either $10$ layers in the whole domain (solid black line) or configurations \eqref{eq:marea_nvar1}-\eqref{eq:marea_nvar4}. Profiles are taken at points $x=4$ km and $x = 16$ km. Dash-dotted red lines are the profiles with constant density.}}
\end{figure}

\begin{figure}[!ht]
	\begin{minipage}{16cm}
		\hspace{-3cm}	\centering\includegraphics[width=1.\textwidth]{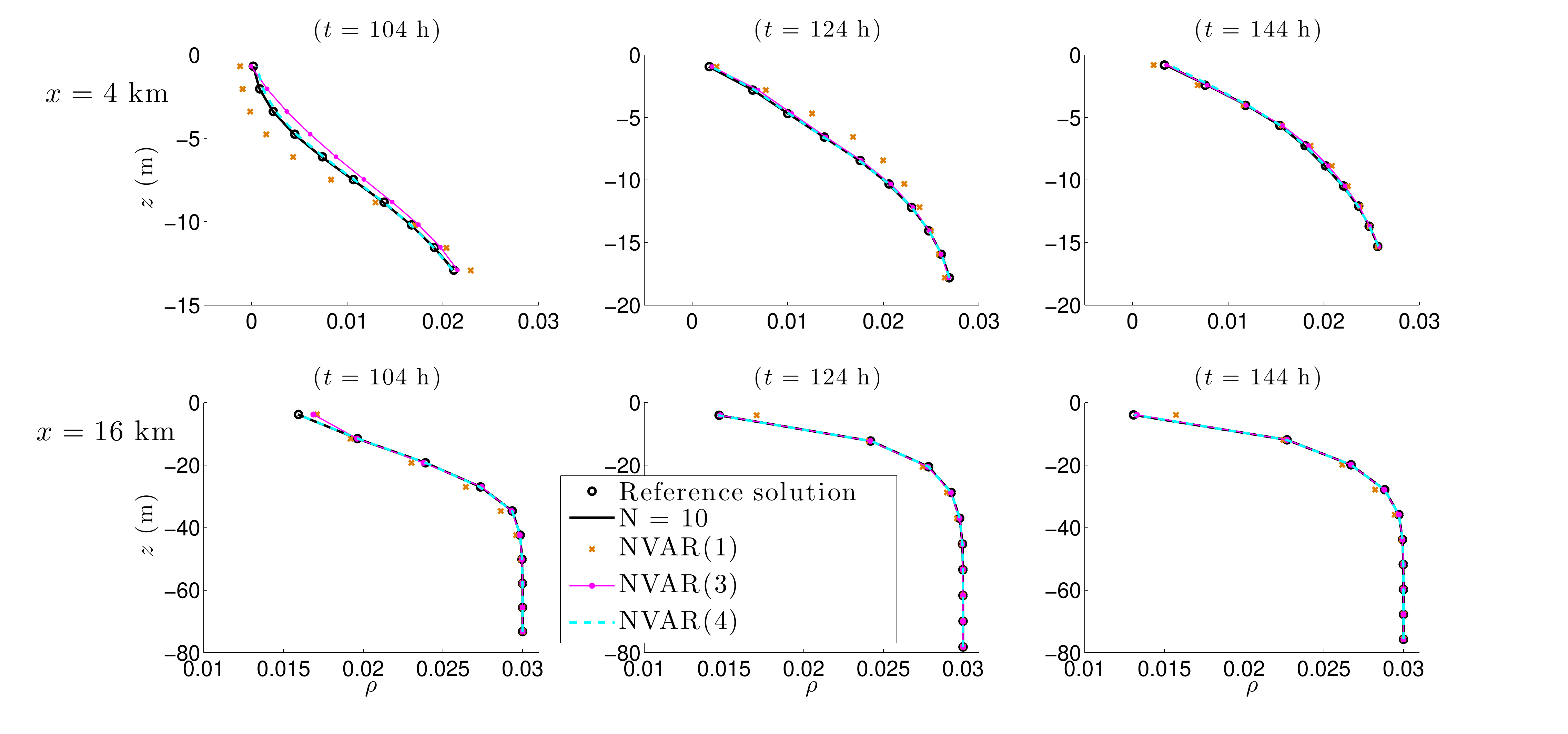}
	\end{minipage}
	\caption{\label{fig:marea_nvar_rho} \footnotesize \textit{Vertical profiles of density obtained with the reference solution (black circles), and the IMEX with $\dt = 10$ s, and either $10$ layers in the whole domain (solid black line) or configurations \eqref{eq:marea_nvar1}-\eqref{eq:marea_nvar4}. Profiles are taken at points $x=4$ km and $x = 16$ km.}}
\end{figure}

Figure \ref{fig:marea_nvar_vel} and \ref{fig:marea_nvar_rho} show the vertical profiles of horizontal velocity and density at points $x=4$ km and $x=16$ km (the top of the peak) at different times. For $x=16$ km, which is far from the zone with reduced number of layers, the vertical profiles for configurations \eqref{eq:marea_nvar3}-\eqref{eq:marea_nvar4} coincides with the ones using a constant number of layers. The approximation with configuration \eqref{eq:marea_nvar1} slightly differs from these ones. For $x=4$ km, which is close to the zone where the transition between the area with constant and variable number of layers occurs, we need configurations \eqref{eq:marea_nvar3}-\eqref{eq:marea_nvar4} to reproduce the profiles with constant number of layers, although the mean variation is well reproduced by all the configurations. As conclusion, we see that configuration \eqref{eq:marea_nvar4} notable reduce the degrees of freedom of the system and it reproduces perfectly the profiles obtained with $10$ layers in the whole domain.

\begin{figure}[!ht]
	\centering\includegraphics[width=0.89\textwidth]{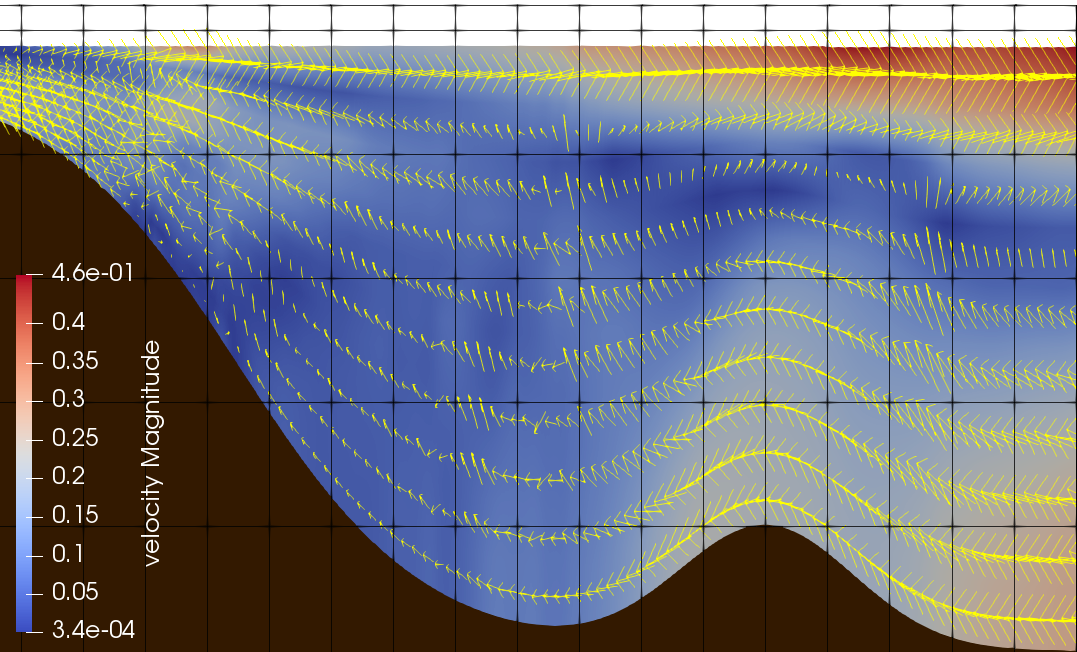}\\[3mm]
	\centering\includegraphics[width=0.89\textwidth]{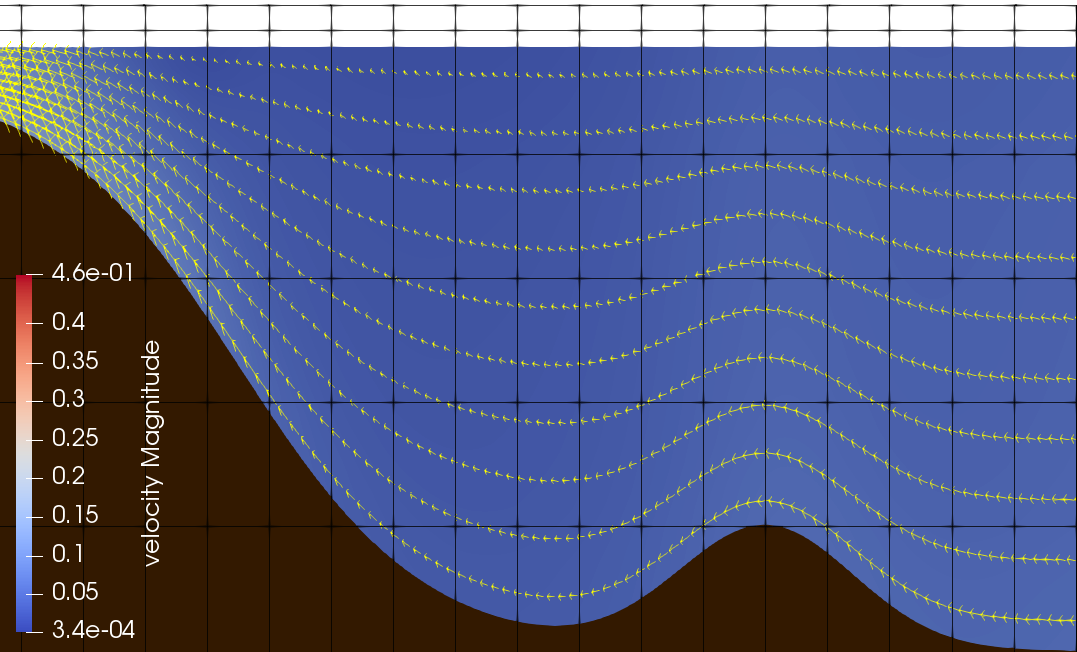}
	\caption{\label{fig:marea_vel_vectors} \footnotesize \textit{Vector map of the velocity field $(u,w)$ at time $t=130$ h for the variable density  (upper figure) and the constant density (lower figure) cases. Colors represent the magnitude of the velocity.}}
\end{figure}

Figure \ref{fig:marea_nvar_vel} also shows the vertical profiles of velocity obtained with constant density, i.e. $\rho_{\a}=0$, in the whole domain. We see that the density field notably changes the profiles of velocity, increasing the vertical structure and dynamics of the flow. 
This  is also observed in Figure 
\ref{fig:marea_vel_vectors}, where the vector velocity field is represented for the variable and constant density configurations. We see again that the vertical structure of the flow increases because of the density field and that the magnitude of the velocity is greater. 

\section{Conclusions}
\label{se:conclu} \indent
The numerical methods proposed in  \cite{bonaventura:2018b} for the barotropic, constant density and hydrostatic case
have been extended to the variable density case in the  Boussinesq regime. To this aim, a transport equation for a variable, which represents the relative deviation of a reference density, is coupled to the mass and momentum equations. 
Although multilayer systems with variable density have been considered previously in the literature, they have never been discretized using a semi-implicit method. An IMEX method is combined  with a specific and consistent discretization of the density equation  and with a multilayer description in which the number of vertical layer can vary along the computational domain.  
An analysis of the linearized multilayer system is presented, showing that the system is hyperbolic for moderate values of the mass transfer terms, while stronger vertical shear induces a loss of hyperbolicity, as can be expected in  a three-dimensional hydrostatic flow. This analysis allows us to define appropriate Courant numbers taking into account the density field.

Some classical tests for variable density flows, as the lock exchange problem, have been performed. We have shown that the proposed semi-implicit method allows us to notably reduce the computational cost of the simulations without a significant loss of accuracy. In particular, we have shown a realistic configuration of variable density flow, which simulates the saltwater intrusion into a river. We have seen that the multilayer configuration can be adapted to complex bathymetries by changing the number of vertical layers if needed, without a loss of accuracy with respect to the simulations with constant number of vertical layers. 

In future work, we will investigate even more flexible and dynamical multilayer discretizations, allowing the number of vertical layers  to vary in time, as well as including  more sophisticated turbulence models accounting for the non constant density field. 

\section*{Acknowledgements}
This work was partially supported by the Spanish Government and FEDER through the research project RTI2018-096064-B-C22. The authors would like to thank Enrique D. Fern\'andez-Nieto and Gladys Narbona-Reina for the interesting discussions related to this work.

\bibliographystyle{plain}
\bibliography{multilayer}

\end{document}